\documentclass[12pt,a4paper]{article}

\usepackage[latin2]{inputenc}
\usepackage[mathscr]{eucal}
\usepackage{t1enc}
\usepackage{amssymb}
\usepackage{amsfonts}
\usepackage{amsmath}

\usepackage{amsthm}
\usepackage{latexsym}
\usepackage[dvips]{color}

\usepackage{graphicx}

\usepackage{bbm}

\newtheorem{Theorem}{Theorem}[section]
\newtheorem{Remark}[Theorem]{Remark}
\newtheorem{Lemma}[Theorem]{Lemma}

\newtheorem{Proposition}[Theorem]{Proposition}
\newtheorem{Conjecture}[Theorem]{Conjecture}
\newtheorem{Corollary}[Theorem]{Corollary}

\begin{document}

\title{\Huge
{\bf Sharp phase transition \\ and critical behaviour in \\ 2D divide and colour models}}

\author{
{\bf Andr\'as B\'alint},
\thanks{E-mail: abalint@few.vu.nl}
{\bf \, Federico Camia},
\thanks{Research supported in part by a VENI grant of the NWO
(Dutch Organization for Scientific Research).}\,
\thanks{E-mail: fede@few.vu.nl}
{\bf \, Ronald Meester}
\thanks{Research supported in part by a VICI grant of the NWO
(Dutch Organization for Scientific Research).}\,
\thanks{E-mail: rmeester@few.vu.nl}\\
{\sl Department of Mathematics, Vrije Universiteit Amsterdam}
}

\maketitle

\begin{abstract}
Consider subcritical Bernoulli bond percolation with fixed parameter $p< p_c$.
We define a dependent site percolation model by the following procedure:
for each bond cluster, we colour all vertices in the cluster black with
probability $r$ and white with probability $1-r$, independently of each other.
On the square lattice, defining the critical probabilities for the site model and its dual,
$r_c(p)$ and $r_c^*(p)$ respectively, as usual, we prove that $r_c(p)+r_c^*(p)=1$ for
all subcritical $p$.
On the triangular lattice, where our method also works, this leads to $r_c(p)=1/2$,
for all subcritical $p$.
On both lattices, we obtain exponential decay of cluster sizes below $r_c(p)$,
divergence of the mean cluster size at $r_c(p)$, and continuity of the percolation
function in $r$ on $[0,1]$.
We also discuss possible extensions of our results, and formulate some natural conjectures.
Our methods rely on duality considerations and on recent extensions of the classical
RSW theorem.
\end{abstract}



\section{Introduction}

\subsection{Definition of the model and main results}

Despite the vast literature on two-dimensional percolation and the tremendous progress made in
its analysis since its introduction as a mathematical theory in \cite{BH}, the exact value of
the critical density is known only for a handful of models.
The latter cases are typically Bernoulli (\textit{independent}) percolation models endowed with
certain \textit{duality} properties, which play a crucial role in the determination of the critical
point.
In this paper, we will be concerned with the study of (the value of) the critical point and the
``phase diagram'' of certain two-dimensional \textit{dependent} percolation models.

Our main object of interest is the two-dimensional \textit{Divide and Color (DaC) model} introduced
by H\"aggstr\"om~\cite{DaC}.
For our purposes, it will be sufficient to consider the simplest version of the model, which can be
described as follows.
Given a graph $G=(\mathcal{V},\mathcal{E})$ with vertex set
$\mathcal{V}$ and edge set $\mathcal{E}$, assign to each edge $e \in \mathcal{E}$ value $1$
(present/open) with probability $p$ and value $0$ (absent/closed) with probability $1-p$, independently
of all other edges.
Denote the resulting $\{0,1\}^{\mathcal{E}}$-valued configuration by $Y$, and the corresponding distribution
by $\nu_p$.
Call \textit{$p$-clusters (``protoclusters'')} the maximal connected components of the graph obtained
by removing from $G$ all the closed edges. Next, colour the vertices of each $p$-cluster black with
probability $r$ and white with probability $1-r$, independently of all other $p$-clusters.
Denote by $X$ the resulting configuration of black and white vertices, by $\mathbb{P}_{p,r}$ the
corresponding distribution, and call \textit{$r$-clusters} the maximal connected
(via the edge set $\mathcal{E}$) white and black subsets of the vertex set $\mathcal{V}$.
These are for us the ``real clusters'', whose percolation properties we are interested in.
Indeed, we will be mainly interested in the percolation properties of $X$ for fixed density $p$ of
open edges, and will consider the edge configuration $Y$ as an auxiliary object, needed to define $X$.
We will in fact argue later in the introduction that our results should still be valid if the
product measure $\nu_p$ on configurations of edges -- corresponding to Bernoulli bond percolation --
is replaced by some other measure with analogous properties of translation invariance and
ergodicity -- e.g., the random cluster measures of which Bernoulli bond percolation is a special
case. Nonetheless, our arguments make substantial use of properties of the product measure
$\nu_p$, and cannot be immediately applied to models defined using other measures.

\begin{figure}[h]
\centering
\includegraphics[scale=0.25, trim= 0mm 0mm 0mm 0mm, clip ]{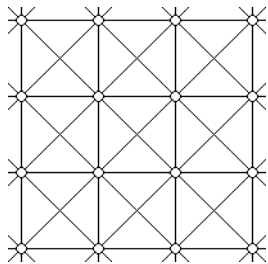}
\caption{Part of the matching graph of the square lattice.}
\label{duallat}
\end{figure}

We first restrict our attention to the square lattice, with vertex set ${\mathbb Z}^2$ and edge set
$\mathcal{E}^2$
given by the edges between nearest neighbour elements of ${\mathbb Z}^2$, and its matching graph,
with vertex set ${\mathbb Z}^2$ and edge set given by the edges between nearest and next-nearest
neighbour elements of ${\mathbb Z}^2$ (i.e., the previous graph with two edges added to each
face of the graph along the two diagonals -- see Figure \ref{duallat}).
We remind the reader that the measure $\nu_p$ on the square lattice has a percolation phase
transition at $p=1/2$~\cite{Kesten}.
We denote by $\Theta(p,r)$ the probability that in the DaC model with parameters $p$ and $r$
the origin $(0,0)$ of the square lattice is contained in an infinite black $r$-cluster, and by
$\Theta^*(p,1-r)$ the probability that it is contained in an infinite white $*$-cluster, where
a $*$-cluster is a connected set of vertices of the matching graph of the square lattice (i.e.,
connections along the diagonals are allowed).
For fixed $p$, we let $r_c(p)=\sup{\{r:\Theta(p,r)=0\}}$ and $r^*_c(p)=\sup{\{r:\Theta^*(p,r)=0\}}$.
In \cite{DaC}, Theorem 2.6, it is shown that $r_c(p)$ and $r_c^*(p)$ are non-trivial.
For fixed $p$ and $r$, we call the model \textit{critical} if (i) $\Theta(p,r)=0$ and
(ii) the mean size of the black $r$-cluster of the origin is divergent (we call the size of a
cluster $C$ its cardinality $|C|$, i.e.\ the number of vertices in the cluster).
In this context, we have the following results.

\begin{Theorem} (Duality) \label{duality}
For all $p<1/2$, $r_c(p)+r_c^*(p)=1$.
\end{Theorem}

We remark that van den Berg has recently proved~\cite{Robnew} that this relation holds
for a large class of percolation models using methods different from those of this paper.
However, the DaC model does not seem to fit in the framework treated in~\cite{Robnew}.

\begin{Theorem} (Exponential decay) \label{exponential-decay}
For $p<1/2$ and $r<r_c(p)$, the size of the black $r$-cluster $C_0^r$ of the origin has an exponentially
decaying tail, i.e. there exists a constant $c(p,r)>0$ such that
\begin{displaymath}
\mathbb{P} _{p,r}(|C_0^r|\geq n)\leq e^{-c(p,r)n}
\end{displaymath}
for all $n\in \mathbb{N}$.
\end{Theorem}
The proof of this result is quite similar to the proof of Theorem 2 in \cite{Voronoi} (see also the
proof of Theorem 5 in \cite{shpthr}), and we do not give it here.

\begin{Theorem} (Criticality) \label{criticality}
The DaC model is critical for $0 \leq p < 1/2$ and $r=r_c(p)$, and for $p=1/2$ and $r \in (0,1)$.
It is not critical for $1/2 < p \leq 1$, where $\Theta(p,r)>0$ for all $r>0$.
\end{Theorem}

Theorem~\ref{duality} amounts to a duality relation between black percolation and white $*$-percolation.
Together with the other two theorems, it provides a complete picture of the phase diagram of the DaC
model, summarized below.

\begin{Corollary} (Phase diagram) \label{phase-diagram}
\begin{itemize}
\item For all $p<1/2$, there exists $r_c(p) \in [1/2,1)$ such that:
\begin{enumerate}
\item If $r<r_c(p)$, there exists an infinite white $*$-cluster a.s.\ and the size of
the black $r$-cluster of the origin has an exponentially decaying tail.
\item If $r=r_c(p)$, $\Theta(p,r_c(p))=\Theta^*(p,1-r_c(p))=0$ and the mean size of the black
$r$-cluster of the origin and of the white $*$-cluster of the origin are infinite.
\item If $r>r_c(p)$, there exists an infinite black $r$-cluster a.s.\ and the size of the
white $*$-cluster of the origin has an exponentially decaying tail.
\end{enumerate}
\item For $p=1/2$, $\Theta(1/2,r)=0$ for all $r \in (0,1)$ and the mean size of the
black $r$-cluster of the origin is infinite.
\item For all $p>1/2$, $\Theta(p,r)>0$ for all $r \in (0,1]$.
\end{itemize}
\end{Corollary}

It is interesting to notice that the two regions of the phase diagram where $\Theta(p,r)>0$,
namely (1) $0 \leq p < 1/2$ and $r>r_c(p)$, and (2) $1/2 < p \leq 1$ and $r>0$, have different properties.
In the first one, there is an infinite black $r$-cluster (somewhere) with probability $1$.
This follows from the ergodicity of the measure $\mathbb{P}_{p,r}$ with respect to translations
when $p<1/2$ (see Section~\ref{prelim}), and the fact that the event that $x \in {\mathbb Z}^2$
belongs to an infinite black $r$-cluster is translation invariant, and has strictly positive probability.
In the second region, where $p>1/2$, the probability
that there is an infinite black $r$-cluster is bounded away from $1$ for all $r<1$.
This also shows that when $p>1/2, r\in(0,1)$, the measure $\mathbb{P}_{p,r}$ is
not ergodic with respect to translations (because of the presence of a unique infinite $p$-cluster).

Another result that follows easily from Theorem~\ref{criticality} is the continuity of the
percolation function $\Theta(p,r)$ as a function of $r$ for $p<1/2$.

\begin{Corollary} \label{continuity}
For all $p<1/2$, $\Theta(p,r)$ is a continuous function of $r$.
\end{Corollary}

The methods used to prove our main results above are not restricted to the square lattice.
In particular, they can be applied to the DaC model on the triangular lattice to obtain the
following theorem, where $p_c(\mathbb{T}):=2\sin(\pi/18)$ is the critical density for Bernoulli bond
percolation
on the triangular lattice~\cite{wierman} (see also \cite{se1,se2}).

\begin{Theorem} (Critical point) \label{critical-point}
In the context of the DaC model on the triangular lattice, for all $p<p_c(\mathbb{T})$, $r_c(p)=1/2$.
\end{Theorem}

We remark that Theorem~\ref{critical-point} looks stronger than Theorem~\ref{duality} because
of the self-duality of site percolation on the triangular lattice, with $*$-clusters being of
the same nature as $r$-clusters, which immediately implies $r_c(p)=r^*_c(p)$.
The proofs of Theorems~\ref{exponential-decay} and~\ref{criticality} and Corollaries~\ref{phase-diagram}
and~\ref{continuity} can also be easily adapted to the triangular lattice.

The DaC model discussed in this section is one of the simplest dependent percolation models
that one can think of. Indeed, it is defined using only product measures. Nonetheless, despite
its simplicity, even results that are by now considered standard for Bernoulli percolation, such
as those discussed in this paper, appear to be much harder to prove for the DaC model than in the
independent case. The proofs of such results, although based on ideas developed for Bernoulli
percolation, require various original arguments that could potentially be helpful in analysing
other dependent percolation models.

The structure of the DaC model is very similar in spirit to that of the random cluster model.
Following this analogy,
the DaC model can be generalised and seen as a particular member of a larger family of models, as
explained in the next section.

\subsection{Other models}

As mentioned just before Theorem~\ref{critical-point}, our methods are robust in the sense that
they work, with obvious modifications, on different lattices.
Another natural extension of our results would be to replace the product measure $\nu_p$ with
other measures.
In particular, we have in mind the class of random cluster measures (of which $\nu_p$ is a special
case -- see, e.g., \cite{grimmett2}).
These are dependent percolation models that unify in a single two-parameter family a variety of
stochastic processes of significant importance for probability and statistical physics, including
Bernoulli percolation, Ising and Potts models.
They are characterized by two parameters, $0 \leq p \leq 1$ and $q>0$, with $q=1$ corresponding
to the Bernoulli percolation measure $\nu_p$.
For $q \geq 1$, they have positive correlation and are believed to have a phase transition in $p$
with exponential decay of correlations for $p<p_c(q)$.
For a detailed account on the random cluster model, the reader is referred to~\cite{grimmett2}.

\begin{Conjecture} \label{random-cluster-measures}
Let $\tilde X = \tilde X(p,q,r)$ be a configuration generated as in the DaC model but with $\nu_p$
replaced by a random cluster measure with parameters $q>0$ and $p<p_c(q)$.
For this model, we conjecture that $r_c(p)+r^*_c(p)=1$ on the square lattice and $r_c(p)=1/2$ on the
triangular lattice.
Furthermore, we conjecture that results analogous to Theorems~\ref{exponential-decay} and~\ref{criticality}
and Corollaries~\ref{phase-diagram} and~\ref{continuity} also hold.
\end{Conjecture}

One major obstacle in proving the conjecture using the methods of this paper is the lack of a proof
of exponential decay of correlations for random cluster measures with $q \neq 1,2$.
The case $q=1$ is treated in this paper. The case $q=2$ is particularly interesting since for $r=1/2$
it corresponds to the Ising model; we address this case in a forthcoming paper.
On the triangular lattice, where self-duality holds, it is easy to see that, at least for $q\geq 2$, the
model defined by $\tilde X$ is critical, in our definition of the term, at the purposed critical point.

\begin{Proposition} \label{universality}
Let $\tilde X(p,q,r)$ be the model defined above on the triangular lattice with $q\geq 2$, $p<p_c(q)$, and
$r=1/2$.
Then there is no infinite $r$-cluster a.s.\ but the mean size of the black $r$-cluster of the origin is
infinite.
\end{Proposition}

Proposition~\ref{universality} implies that $r_c(p) \geq 1/2$ and suggests that indeed $r_c(p)=1/2$.
It has a ``universality'' flavour since it suggests that self-duality alone (almost) determines the
critical value $r_c(p)=1/2$.
One should however compare this result to Theorem~\ref{criticality}, which shows that when $p=p_c$
(i.e., either $1/2$ on the square lattice or $2\sin(\pi/18)$ on the triangular lattice) the DaC model
has a critical segment, $r \in (0,1)$, rather than a single critical point.
The difference between the two cases is that when $p<p_c$ the size of $p$-clusters has an exponentially
decaying tail, while at $p_c$ the mean $p$-cluster size diverges and $p$-clusters form circuits around the
origin at all scales, turning the percolation of $r$-clusters effectively into a one-dimensional problem.

\subsection{Scaling limits}

It is natural to ask what the continuum scaling limit (when the lattice spacing is sent to zero)
of the DaC model is on the two critical ``curves''
(1) $p \in [0,p_c)$, $r=r_c(p)$ and (2) $p=p_c$, $r \in (0,1)$.

For the first critical curve, based on universality considerations, we expect the scaling limit
to be the same as for critical Bernoulli site percolation, corresponding to $p=0$ and $r=r_c(0)$.
In particular, we expect crossing probabilities to converge to Cardy's formula~\cite{cardy} (as
proved by Smirnov~\cite{smirnov} for critical Bernoulli site percolation on the triangular lattice)
and the set of all interfaces between black $r$-clusters and white $*$-clusters to converge to the
Continuum Nonsimple Loop process described in~\cite{cn1,cn2}.
This is in line with the general principle that short range correlations, as produced by the
$p$-clusters below $p_c$, do not affect the critical behaviour and the scaling limit.
The DaC model in this regime can be seen as Bernoulli site percolation on a random graph whose
vertices are the $p$-clusters, and as long as $p<p_c$, the random graph will be, in some sense,
``close'' to the underlying regular lattice.
In other words, under the action of the renormalisation group, the critical curve (1) should
have a unique fixed point, namely, $p=0$, $r=r_c(0)$.

We expect similar considerations to hold when the product measure $\nu_p$ is replaced by a
different random cluster measure (or even in greater generality), and make the following natural
conjecture, stated for simplicity for the triangular lattice.

\begin{Conjecture}
Let $\tilde X(p,q,r)$ be the model of Proposition~\ref{universality}. For all $q>0$ and all
$p<p_c(q)$, the (site percolation) scaling limit of $\tilde X(p,q,1/2)$ is the same as the
scaling limit of critical Bernoulli (site) percolation.
\end{Conjecture}

In the case of the second critical curve, we expect a different situation,
with different scaling limits for different values of $r$.
We expect, for instance, that the scaling limit of crossing probabilities will depend on $r$
and will not in general be given by Cardy's formula.

\subsection{Strategy of the proof of Theorem~\ref{duality}}

The proof of our main result, Theorem~\ref{duality}, follows a ``modern'' version
(using Russo's formula) of the celebrated proof~\cite{Kesten} by Kesten that the critical
probability for Bernoulli (independent) bond percolation on the square lattice is $1/2$
(see also~\cite{Russocrpr}).
However, since we are dealing with a \textit{dependent} percolation model, the proof requires
various modifications, needed for instance to avoid gathering ``too much information.''
Also, in this modified version of Kesten's strategy, and due to the dependence
structure of the DaC model, we cannot apply the ``traditional'' RSW theorem.
We will instead use a recent version of it taken from~\cite{BBB}, which is a strengthened
form of the RSW type theorem in~\cite{Voronoi}.

We now describe briefly (and somewhat imprecisely) what one would do in the case of Bernoulli
percolation, corresponding to $p=0$. Some notational remarks first:
we shall omit the subscript $r$ from the notation of the
DaC measure $\mathbb{P}_{p,r}$ and also often write $r_c$ and $r_c^*$ for $r_c(p)$ and $r_c^*(p)$
respectively if no confusion is possible.

It follows from standard arguments that $1-r_c^*\leq r_c$.
Therefore, the main task is to prove that this inequality is not strict.
We proceed by contradiction, assuming that the open interval $(1-r_c^*,r_c)$ is non-empty.
Let $V^b_{n,3n}$ denote the presence of a vertical black crossing of an $n \times 3n$
rectangle and $H^{w*}_{n,3n}$ the presence of a horizontal white $*$-crossing of the same
rectangle (precise definitions will be given in Section~\ref{basic-defs} below).
Accordingly, let $H^b_{n,3n}$ denote the event that there is a horizontal black crossing
of an $n \times 3n$ rectangle, $V^{w*}_{n,3n}$ that there is a vertical white $*$-crossing
of the same rectangle.
Since there is a.s.\ no percolation of black vertices for $r \in (1-r^*_c,r_c)$,
it is easy to prove that $\limsup_{n \to \infty} {\mathbb P}_r (V^b_{n,3n}) < 1$.
But then, by the simple but crucial observation that, no matter how one chooses to colour
the vertices inside the $n \times 3n$ rectangle, there is always either a vertical black
r-crossing or a horizontal white $*$-crossing, $\liminf_{n \to \infty}{\mathbb P}_r(H^{w*}_{n,3n})>0$
for all $r \in (1-r_c^*,r_c)$.
This implies that there is a uniform positive probability to have a horizontal white $*$-crossing in
the lowest half of an $n \times 8n$ rectangle for all $n$ large enough. (We note that Kesten
used squares in his proof. In our case, due to the dependence structure, it will be more
convenient to use rectangles. To avoid confusion and prepare the reader for what will
come, we employ the same rectangles here.)

Consider the lowest such crossing $\pi_1$ and look for a black vertical $r$-crossing in
the left half of the same $n \times 8n$ rectangle from the top of the rectangle to $\pi_1$.
Since there is no white $*$-percolation for $r \in (1-r^*_c,r_c)$, with the help of the RSW
theorem, one can show that
a black vertical crossing from the top of the rectangle to
$\pi_1$ exists with probability bounded away from zero.
Consider the leftmost such crossing $\pi_2$.
Due to the properties of lowest crossings, the presence of such a black crossing
implies the presence of a white vertex $x$ on $\pi_1$ which is pivotal for the event
$H^{w*}_{n,8n}$.

Next, center at $x$ a sequence of nested annuli intersecting the $n \times 8n$ rectangle.
Inside the portion of each annulus intersecting the rectangle and lying above $\pi_1$ and
to the right of $\pi_2$, look for a black crossing joining $\pi_2$ with $\pi_1$.
Once again, the existence of such crossings with uniform positive probability is assured by the
RSW theorem.
Every such crossing gives another white vertex on $\pi_1$ which is pivotal for $H^{w*}_{n,8n}$.
In this way, choosing the annuli appropriately, one can find many pivotal vertices with high
probability.
Using Russo's formula it is then possible to conclude that  ${\mathbb P}_r (H^{w*}_{n,8n})$ has
a very large (negative) derivative for all $r \in (1-r_c^*,r_c)$, obtaining a contradiction.

The argument above relies on properties of lowest and leftmost crossings, and in
particular uses the fact that a lowest (respectively, leftmost) crossing can be
found without exploring the area above (resp., to the right of) the crossing itself.
In the case of Bernoulli percolation this implies that the configuration above the lowest
crossing can be coupled to an independent configuration, and the probability to find a black
r-crossing in the left half of the rectangle can be bounded below using the RSW theorem.
The same type of argument applies to the portions of annuli to the right of $\pi_2$ and
above $\pi_1$, where the configurations can again be coupled to independent configurations,
and the probabilities of finding the appropriate crossings bounded below once again
using the RSW theorem.

In our case, similar arguments can be used, but the dependence in the model makes them
significantly more complex.
Moreover, as remarked above, the ``traditional'' RSW theorem cannot be used, and we have
to resort to a more recent version~\cite{BBB} which is weaker but more general and, as it
turns out, still sufficiently strong for our purposes.

To deal with the dependence structure of the model, in some situations we will ``fatten''
certain collections of vertices (e.g., vertices forming a crossing) by adding to them their
$p$-clusters. This procedure identifies closed ``barriers'' of edges with the property that
colour configurations on different sides of a barrier are conditionally independent
(conditioned on the barrier).

We will also use algorithmic constructions carefully designed to explore certain domains
looking for monochromatic crossings without obtaining too much information. This will allow
us to couple in a useful way the DaC measure ${\mathbb P}_{p,r}$ conditioned on some specific
$\sigma$-algebras corresponding to the information obtained while looking for crossings with
an unconditional version.

\subsection{Outline of the paper}

In Section \ref{prelim}, we present the definitions and introduce notation.
Then, we collect the tools which are needed to prove the main results.
These include known results such as the exponential decay property of subcritical
Bernoulli bond percolation,
the FKG inequality for the measure $\mathbb{P} _{p,r}$,
and the modern RSW theorem from \cite{BBB}.
Then we give the natural analogue of Russo's formula for the DaC model
(Theorem \ref{tRusso}), and finally state that
percolation occurs with positive probability if and only if certain
rectangles can be crossed with high probability (Lemma \ref{lessoreq}).

Section \ref{mainproof} contains the proof of Theorem \ref{duality}.
In Section \ref{consequences}, we prove Theorem \ref{criticality} and
Corollaries \ref{phase-diagram} and \ref{continuity}.
Finally, in Section \ref{triang}, we prove our results on the triangular lattice.

\section{Preliminaries} \label{prelim}

\subsection{Basic definitions and notation} \label{basic-defs}

We consider the square lattice, with vertices the points of $\mathbb{Z}^2$, and
edges between adjacent vertices (that is, between vertices at Euclidean distance 1).
With the usual abuse of notation, we denote both the graph and its vertex set by $\mathbb{Z}^2$,
and we write $\mathcal{E}^2$ for the edge set of this graph.

The state space of our configurations is defined as $\Omega :=\Omega _{D}\times \Omega _{C}$,
where $\Omega _{D}:=\{0,1\}^{\mathcal{E}^2 }$ corresponds to Bernoulli bond percolation,
and $\Omega _{C} :=\{0,1\}^{\mathbb{Z}^2}$ corresponds to colouring.
We identify 0 with the colour white, and 1 with black.
The probability measure $\mathbb{P}_{p,r}$ is the measure (on the usual $\sigma$-algebra on $\Omega$)
obtained by the procedure described in the introduction.

We introduce the set $\tilde\Omega \subset \Omega$ as the set of configurations such that vertices in the same
$p$-cluster have the same colour, and we equip $\tilde\Omega $ with a partial ordering as follows.
For $\omega _1=(\eta _1,\xi _1), \omega _2=(\eta _2,\xi _2) \in \tilde\Omega $ we say that
$\omega _1 \geq \omega _2$ if, for all $x\in \mathbb{Z}^2$, we have $\xi _1 (x)\geq \xi _2(x)$.
Note that the ordering depends on the colours of the vertices only, not on the bond configurations. All the
configurations in this paper are silently assumed to be in $\tilde{\Omega}$.
We call an event $A \subset \tilde{\Omega}$ \textit{increasing} if $\omega \in A$ and
$\omega '\geq \omega $ implies $\omega '\in A$.
$A$ is a \textit{decreasing event} if $A^c$ is increasing.

We call a sequence of vertices $(x_0,x_1,\ldots ,x_n)$ in $\mathbb{Z}^2$ a \textit{(self-avoiding) path}
if for all $i=0,\ldots ,n-1$, $x_i$ and $x_{i+1}$ are adjacent, and for any
$0\leq i< j\leq n,\hspace{0.1cm} x_i\neq x_j$.
The definition of a \textit{$*$-path} is similar, but for $i=0,\ldots ,n-1$,
the vertices $x_i$ and $x_{i+1}$ need to be just \textit{$*$-adjacent} instead of adjacent,
which means that their Euclidean distance is 1 or $\sqrt{2}$.
A ($*$-)\textit{circuit} is defined in the same way as a ($*$-)path except that $x_n=x_0$.
A \textit{horizontal crossing} of a rectangle $R=[a,b]\times [c,d]$, with $a,b,c,d \in {\mathbb Z}$,
is a path ${x_0,x_1,\ldots ,x_n}$ such that $x_0\in \{a\}\times [c,d]$, $x_n\in \{b\}\times [c,d]$
and for all $i$, $x_i\in R$.
A \textit{vertical crossing} of the same rectangle is a path ${x_0,x_1,\ldots ,x_n}$ such
that $x_0\in [a,b]\times \{ d\}$, $x_n\in [a,b]\times \{c\}$ and for all $i$, $x_i\in R$.
\textit{$*$-circuits}, \textit{horizontal $*$-crossings}, and \textit{vertical $*$-crossings} are defined
by replacing paths by $*$-paths in the above definitions.

A \textit{black path} is a path $\pi ={x_0,x_1,\ldots ,x_n}$ such that for all $i=0,\ldots ,n,$ $x_i$
is black (i.e. $X(x_i)=1$).
\textit{Black circuits}, \textit{black horizontal crossings}, \textit{black vertical crossings} are defined
analogously.
A \textit{black cluster} is a maximal subset $K$ of $\mathbb{Z}^2$ such that between any two vertices of $K$
there exists a black path.
The definitions of \textit{white path}, \textit{white circuit}, \textit{white horizontal crossing},
\textit{white vertical crossing}, \textit{white cluster} are obtained by replacing black with white.
Black and white ($*$-)paths, ($*$-)circuits, ($*$-)crossings, and ($*$-)clusters are defined analogously.

Let $S_{n,m}$ denote the rectangle $[0,n]\times [0,m]$, with $n,m \in {\mathbb N}$.
Denote by $V^b_{n,m}$ the event that there is a vertical black crossing in the rectangle $S_{n,m}$;
let $H_{n,m}^b$ be the corresponding event with a horizontal crossing.
Furthermore, let $B_n^b$ denote the event that there is a black circuit surrounding the midpoint in the
annulus $A_n :=S_{3n,3n} \setminus \left( S_{n,n}+\left( n,n\right) \right)$.
Here and later, for a set $S$ and a vector $v$, we use the notation $S+v:=\{x:x-v\in S\}$.
The analogous events with white crossings are denoted by $V_{n,m}^w$, $H_{n,m}^w$, and $B_n^w$,
respectively.
A $*$ in the notation will indicate that we are referring to $*$-crossings and $*$-circuits -- for example,
$V^{w*}_{n,m}$ denotes the event that there is a vertical white $*$-crossing in $S_{n,m}$.

Let $d$ denote the $L_1$ distance.
The distance between two sets of vertices $V_1$ and $V_2$ is defined by $d(V_1,V_2):=\min
\{d(x,y)\hspace{0.05cm}:\hspace{0.05cm x\in V_1,y\in V_2}\}$.
Let $\partial B(v,n)$ denote the circle of radius $n$ with center at vertex $v$ in the metric $d$,
i.e., $\partial B(v,n):=\{w:d(v,w)=n\}$.
For a vertex $v\in \mathbb{Z}^2 $, let $C_v ^p$ be the open $p$-cluster of $v$, i.e., the set of
vertices that can be reached from $v$ through edges that are open in the underlying Bernoulli bond
percolation with parameter $p$.
Let us define the \textit{dependence range} of a vertex $v$ by
$\mathcal{D}(v):=\max \{n\in \mathbb{N}:C_v ^p\cap \partial B(v,n) \neq \emptyset\}$.

We call an edge set $E=\{e_1,e_2,\ldots ,e_k\}$ a \textit{barrier} if
removing $e_1,e_2,\ldots ,e_k$ (but not their end-vertices) separates the graph $\mathbb{Z}^2$ into two
or more disjoint connected subgraphs, of which exactly one is infinite.
We call the infinite component of $({\mathbb Z}^2,\mathcal{E}^2) \setminus E$
the \textit{exterior} of $E$, and denote it by $ext(E)$.
We call the union of the finite components the \textit{interior} of $E$, and denote it by $int(E)$.
(Note that a barrier as defined above corresponds to a dual circuit in bond percolation. However, since
we work with a different sort of duality throughout this paper, we adopt a different term to avoid
confusion.)
$E=\{e_1,e_2,\ldots ,e_k\}$ is a \textit{closed barrier} if $E$ is a barrier and $e_i$ is closed
in the Bernoulli bond percolation (i.e. $Y(e_i)=0$).
For a vertex set $A\subset \mathbb{Z}^2$, let $\Delta A$ denote the \textit{edge boundary} of $A$, that is,
$\Delta A:=\{(x,y)\in \mathcal{E}^2 : x\in A, y\in \mathbb{Z}^2 \setminus A\}$.
Note that for $p < p_c$, the edge boundary of any $p$-cluster is a closed barrier.

\subsection{Preliminary results}
\label{prelim-results}

In this subsection we collect several results that are mostly known, follow directly from known results,
or can be proved using variations of classical arguments. The exception is Lemma \ref{stochdom}, which
is new and very important in the forthcoming construction.

\begin{Theorem}\emph{(\cite{menshikov},\cite{ab})}
\label{lexpdecay}
If $p<1/2$, there exists $\psi (p)>0$ such that for all n we have
\begin{displaymath}
\nu _{p}(\mathcal{D}(0)\geq n)< e^{-n\psi (p)}.
\end{displaymath}
\end{Theorem}

A short and simple proof of the following classical result is given in Section 4.2 of \cite{shpthr}.

\begin{Lemma}\label{selfdl}
The restriction to a rectangle $R$ of any colour configuration $\xi \in \Omega_C$ contains
either a black vertical crossing or a white horizontal $*$-crossing of $R$, but never both.
In particular, for any $n,m$ we have
\begin{displaymath}
V_{n,m}^b = \left( H_{n,m}^{w*} \right) ^c.
\end{displaymath}
\end{Lemma}

The proof of the next result, which shows positive correlation for the DaC model, was obtained
by H\"aggstr\"om and Schramm and included in~\cite{DaC}.

\begin{Theorem}\emph{(\cite{DaC})}
\label{FKG}
Let $A,B$ be increasing events. Then, for any $p,r\in [0,1]$,
\begin{displaymath}
\mathbb{P}_{p,r}(A\cap B)\geq \mathbb{P}_{p,r}(A)\mathbb{P}_{p,r}(B).
\end{displaymath}
\end{Theorem}

We shall also need a result which makes precise (and generalizes) the observation that an edge
between two vertices of the same colour is more likely to be open than an edge between vertices
whose colours are unknown.
Let us consider the following scenario for Lemma \ref{stochdom} below:
let $B_1, B_2,\ldots ,B_k$ be barriers,
$e_1,\ldots ,e_l$ edges in $U:=\bigcup _{i=1} ^kint(B_i)$, $u_1,\ldots ,u_m$ vertices in $U$
(where $k,l,m\in \{0,1,2,\ldots \}$), and denote $\bigcap _{i=1} ^kext(B_i)$ by $V$.
Fix states $s_i\in \{0,1\}$, $i=1,\ldots ,l$ and colours $c_j\in \{0,1\}$, $j=1,\ldots ,m$.
Let $v_1,\ldots ,v_n$ be vertices in $V$ (where $n\in \{0,1,2,\ldots \}$), and let
$\kappa \in\{0,1\}$ be a colour.
Let $I$ denote the event that $B_1, B_2,\ldots ,B_k$ are closed,
$Y(e_i)=s_i$, $X(u_j)=c_j$, and $v_1,\ldots ,v_n$ all have colour $\kappa $.

\begin{Lemma} \label{stochdom}
The conditional distribution of the edges in $V$, conditioned on the event $I$ described
above, stochastically dominates the measure $\nu _p$.
\end{Lemma}

\medskip\noindent
{\bf Proof.}
We shall prove the lemma with an iteration, determining the states of edges in $V$ one
after another.
Let $f_1,f_2,\ldots, f_g$ be edges in $V$, and $t_1,t_2,\ldots ,t_g\in \{0,1\}$ states,
where $g\in \{0,1,2,\ldots \}$.
Let us consider the events $A:=\bigcap_{i=1} ^g\{Y(f_i)=t_i\}$, and $J:=I\cap A$.
Take an edge $e$ in $V$ whose state is not determined by $A$.
Note that there is no further restriction on the location of $e$: it may be incident on 0, 1 or 2
vertices from $\{ v_1,\ldots ,v_n \}$. We shall first show that
\begin{equation} \label{eopenJ}
\mathbb{P}_{p,r}(\textrm{$e$ open}|J)\geq p.
\end{equation}

It is easy to see that (\ref{eopenJ}) is equivalent to
$$(1-p)\mathbb{P}_{p,r}(\textrm{$e$ open} , J)\geq p\mathbb{P}_{p,r}(\textrm{$e$ closed} , J),$$
which is also equivalent to
$$(1-p)\mathbb{P}_{p,r}(J|\textrm{$e$ open})\nu _p(\textrm{$e$ open})\geq p\mathbb{P}_{p,r}(J|\textrm{$e$
closed})\nu _p(\textrm{$e$ closed}).$$
Since $\nu _p(e$ open$)=p$ and $\nu _p(e$ closed$)=1-p$, it remains to show that
\begin{equation}\label{givenstateeJ}
\mathbb{P}_{p,r}(J|\textrm{$e$ open})\geq \mathbb{P}_{p,r}(J|\textrm{$e$ closed}).
\end{equation}
This may be seen as follows. Since the first step in constructing a configuration corresponds to
Bernoulli bond percolation, the states of edges other than $e$
are independent of the state of $e$. For $J$ to occur, the edges in $\bigcup _{i=1}^kB_i$ need to be closed.
In that case, the colouring of $u_1,\ldots ,u_m$ is not influenced by the state of $e$, since every vertex
in question is in the interior
of one of the closed barriers.
Therefore, the only thing left to prove is that the probability that the vertices
$v_1,\ldots ,v_m$ all have colour $\kappa $ is greater given that $e$ is open
than given that $e$ is closed.
This follows immediately from a very simple coupling between $\mathbb{P}_{p,r}(J|\textrm{$e$ open})$
and $\mathbb{P}_{p,r}(J|\textrm{$e$ closed})$ in which all the edges except $e$ are in the same
state, since when $e$ is open the number of $p$-clusters that need to be assigned colour
$\kappa $ is smaller than or equal to the number of $p$-clusters that need to be assigned
colour $\kappa $ when $e$ is closed.
This observation proves (\ref{givenstateeJ}), finishing the proof of (\ref{eopenJ}).

The full stochastic domination can be shown using (\ref{eopenJ}) iteratively as follows.
We shall condition on $I$.
Fix a deterministic ordering of the edges in $V$ and use the following iteration:
\begin{enumerate}
\item Start with $A=\{ 0,1 \}^{\mathcal{E}(V)}$ where $\mathcal {E}(V)$ is the set of edges contained in $V$.
\item Determine the state $s_e$ of the first edge $e \in \mathcal {E}(V)$ in the ordering whose state is not
yet determined by $A$, according to the conditional distribution $\mathbb{P}_{p,r}( \, \cdot \, | \, I\cap A)$.
\item $A:=A\cap \{Y(e)=s_e\}$.
\item Go back to step 2.
\end{enumerate}
It is clear that every edge in $V$ gets a state drawn from the correct
distribution after finitely many steps.
On the other hand, we know from (\ref{eopenJ}) that for all $A$, the marginal
of $\mathbb{P}_{p,r}( \, \cdot \, | \, I\cap A)$ on $\Omega_D$ dominates
$\nu_p( \, \cdot \, )$. This proves the desired stochastic domination.
\hfill $\Box $\\

\medskip
In the proof of Theorem~\ref{duality}, we will use an RSW type theorem that was recently obtained
by van den Berg, Brouwer and V\'agv\"olgyi~\cite{BBB}.
This is a stronger version of the RSW type theorem used by Bollob\'as and Riordan in~\cite{Voronoi}.
Such results are weaker than the classical RSW theorem but more general, and can be applied to models
for which the classical RSW theorem has not been proved.
We remark that in our proof of $r_c(p)+r_c^*(p)\leq 1$ in Section \ref{mainproof} we seem to need
the full strength of the result of van den Berg, Brouwer and V\'agv\"olgyi, as the weaker form
proved by Bollob\'as and Riordan does not seem to suffice for our purposes.
Stated for the DaC model, the result reads as follows.

\begin{Lemma}
\label{RSW}
For any $p<1/2,r\in [0,1]$, we have\\
\begin{displaymath}
\textrm{(a) If }\limsup \limits _{n\to \infty }\mathbb{P}_{p,r}(V_{n,\rho n}^b)>0\textrm{ for some $\rho
>0$},
\end{displaymath}
\begin{displaymath}
\textrm{then }\limsup \limits _{n\to \infty }\mathbb{P}_{p,r}(V_{n,\rho n}^b)>0\textrm{ for all $\rho >0$};
\end{displaymath}
\begin{displaymath}
\textrm{(b) If }\limsup \limits _{n\to \infty }\mathbb{P}_{p,r}(V_{n,\rho n}^{w*})>0\textrm{ for some $\rho
>0$},
\end{displaymath}
\begin{displaymath}
\textrm{then }\limsup \limits _{n\to \infty }\mathbb{P}_{p,r}(V_{n,\rho n}^{w*})>0\textrm{ for all $\rho
>0$}.
\end{displaymath}
\end{Lemma}

\medskip\noindent
{\bf Proof.} First we prove (a).
Following~\cite{BBB}, Section 4.3, and~\cite{shpthr}, Section 5.1, it suffices to check conditions
(1)--(5) below.
For a set $R$ and $\lambda \neq 0$, we write $\lambda R$ for $\{x:x/\lambda \in R\}$. We consider the following
five conditions.
\begin{enumerate}
\item[(1)] For any rectangle $R$, if $h$ and $v$ are a horizontal and a vertical crossing of $R$,
respectively, then $d(h,v)\leq 1$.
\item[(2)] Increasing events are positively correlated.
\item[(3)] The model has the symmetries of $\mathbb{Z}^2$, i.e., is invariant under translations by the
vectors $(1,0)$ and $(0,1)$, reflections through the coordinate axes of ${\mathbb Z}^2$, and rotations of
90 degrees.
\item[(4)] Disjoint regions are asymptotically independent as we ``zoom out'' (the precise formulation
of this condition will be given in Lemma \ref{asind}).
\item[(5)] For any fixed rectangle $R$ there is a constant $C$ such that the length of a horizontal
crossing of $\lambda R$ is bounded from above by $\lambda^C$ if $\lambda $ is large enough.
\end{enumerate}

Condition (1) clearly holds here, since horizontal and vertical black crossings of the same rectangle have
at
least one vertex in common.
Condition (2) is given by Theorem \ref{FKG}.
Condition (3) can be checked easily.
For condition (4), see Lemma \ref{asind} below.
Condition (5) obviously holds, since the model is discrete.

For the proof of (b), the same conditions need to be checked with $*$-crossings instead of crossings
in (1) and (5), and increasing events replaced by decreasing events in (2).
The new first condition still holds, since even though a horizontal $*$-crossing and a vertical one
of the same rectangle do not necessarily have a vertex in common, they are at distance at most 1 from
each other.
The new second condition, namely that decreasing events are positively correlated, is an easy consequence
of Theorem \ref{FKG}, since the complement of a decreasing event is an increasing event.
\hfill $\Box $\\

\medskip
The next lemma immediately implies weak mixing and ergodicity
for the DaC model when $p<p_c$.

\begin{Lemma}
\label{asind}
Let $p<1/2,r\in [0,1]$. Then for disjoint rectangles $R_1$ and $R_2$,
for any $\varepsilon >0$ there exists $\lambda _0>0$ such that for all $\lambda >\lambda _0$,
for any events $A_1, A_2$ defined in terms of the colouring of vertices in $\lambda R_1$ and
$\lambda R_2$ respectively, we have
\begin{displaymath}
|\mathbb{P}_{p,r}(A\cap B)-\mathbb{P}_{p,r}(A)\mathbb{P}_{p,r}(B)|\leq \varepsilon .
\end{displaymath}
\end{Lemma}

\medskip\noindent
{\bf Proof.} Fix $p<p_c$.
Let $R_1$ and $R_2$ be rectangles at distance $k>0$. Take arbitrary events $A_1, A_2$ defined in terms
of the colours of vertices in $R_1$ and $R_2$ respectively.
Let $K=K_{R_1,R_2}$ be the event that $R_1$ and $R_2$ are separated by a closed barrier in the
bond configuration.
If $K$ occurs, then the colours of the vertices in $R_1$ and $R_2$ are conditionally independent.
Therefore, $A_1$ and $A_2$ are conditionally independent, conditioned on $K$.
The law of total probability, together with the previous observation, gives

\begin{eqnarray}
\mathbb{P}_{p,r}(A_1\cap A_2) = \mathbb{P}_{p,r} (A_1|K)\mathbb{P}_{p,r}(A_2|K)\nu _p(K)+\mathbb{P}_{p,r}
(A_1\cap A_2|K^c)\nu _p(K^c) \label{A1A2}
\end{eqnarray}

Since

\begin{equation}\label{A1}
\mathbb{P}_{p,r}(A_1)=\mathbb{P}_{p,r} (A_1|K)\nu _p(K)+\mathbb{P}_{p,r} (A_1|K^c)\nu _p(K^c)
\end{equation}
and
\begin{equation}\label{A2}
\mathbb{P}_{p,r}(A_2)=\mathbb{P}_{p,r} (A_2|K)\nu _p(K)+\mathbb{P}_{p,r} (A_2|K^c)\nu _p(K^c),
\end{equation}
by substituting the right hand sides of equations (\ref{A1A2}), (\ref{A1}) and (\ref{A2}), using the
triangle inequality, we obtain

\begin{displaymath}
|\mathbb{P}_{p,r}(A_1\cap A_2)-\mathbb{P}_{p,r}(A_1)\mathbb{P}_{p,r}(A_2)|
\end{displaymath}
\begin{displaymath}
\leq |\mathbb{P}_{p,r} (A_1|K)\mathbb{P}_{p,r} (A_2|K)||\nu _p(K)-\nu _p(K)^2|+|\nu _p(K^c)||Q|,
\end{displaymath}
where $|Q|\leq 4$ since $Q$ is the sum of four products of probabilities.

We also need to notice that if none of the vertices in $R_1\cup R_2$ has a dependence range of at least
(say) $\frac{k}{3}$ in the initial random
bond configuration Y, then $K$ occurs. Therefore,

\begin{eqnarray*}
 \nu _p(K^c)&\leq &  \nu _p\left( \bigcup _{v\in R_1\cup R_2}\{\mathcal{D}(v)\geq \frac{k}{3}\}\right)  \\
&\leq &  (|R_1|+|R_2|)\nu _p (\mathcal{D}(0)\geq \frac{k}{3})\\
&\leq &(|R_1|+|R_2|)e^{-\psi(p)\frac{k}{3}},
\end{eqnarray*}
according to Theorem \ref{lexpdecay}.

It immediately follows that for the probability of the event $L:=K_{\lambda R_1,\lambda R_2}$, we have
$\nu _p(L^c)\leq (|R_1|+|R_2|)\lambda ^2e^{-\psi(p)\frac{\lambda k}{3}}\to 0$ as
$\lambda \to \infty$.\\
Fix $\varepsilon >0$, and choose $\lambda _0$ so large that for $\lambda >\lambda _0$, we have $\nu
_p(L^c)<\frac{\varepsilon }{8} $.
Take arbitrary events $A_1$ and $A_2$, defined in terms of the colours of vertices in $\lambda R_1$ and
$\lambda R_2$, respectively.
Since $\nu _p(L)\geq 1-\frac{\varepsilon }{8}$ implies $|\nu _p(L)-\nu _p(L)^2|\leq \frac{\varepsilon }{4}$,
we obtain
\begin{eqnarray*}
& &|\mathbb{P}_{p,r}(A_1\cap A_2)-\mathbb{P}_{p,r}(A_1)\mathbb{P}_{p,r}(A_2)|\\
& &\leq |\mathbb{P}_{p,r} (A_1|L)\mathbb{P}_{p,r} (A_2|L)|\cdot |\nu _p(L)-\nu _p(L)^2|+|\nu _p(L^c)|\cdot
4\\
& &\leq 1\cdot \frac{\varepsilon }{4}+\frac{\varepsilon }{8}\cdot 4\\
& & < \varepsilon ,
\end{eqnarray*}
proving the lemma.
\hfill $\Box $\\

\begin{Corollary}
\label{erg}
For $p<1/2$, the measure $\mathbb{P}_{p,r}$ is weakly mixing and therefore ergodic with
respect to translations.
\end{Corollary}

We also need a version of Russo's formula~\cite{Russocrpr} (see
also~\cite{Grimmett}).
Let $A$ be an event, and let $\omega =\left( \eta ,\xi \right) $ be a configuration
from $\tilde\Omega$. Let $C$ be an open $p$-cluster from $\eta$.
We call $C$ \textit{pivotal} for the pair $\left( A,\omega \right) $ if $I_A(\omega )\neq I_A(\omega ')$
where $I_A$ is the indicator function of $A$, $\omega' = (\eta,\xi')$, and $\xi'$ agrees with $\xi$
everywhere except that the colour of the vertices in $C$ is different.

\begin{Theorem}
\label{tRusso}
Let $W$ be a set of vertices with $|W|<\infty $, and let $A$ be
an increasing event that depends only on the colours of vertices in $W$.
Then we have, for any $p\in [0,1]$,
\begin{displaymath}
\frac{d}{dr}\mathbb{P}_{p,r}(A)=\mathbb{E}_{p,r}(n(A)),
\end{displaymath}
where $n(A)$ is the number of $p$-clusters which are pivotal for $A$.
\end{Theorem}

\medskip\noindent
{\bf Proof sketch.} Let us denote the (finite) set of partitions of the vertices in $W$
which are compatible with a bond configuration by $\mathcal{P}_W$, and the (random)
partitioning determined by the initial bond percolation by $\mathcal{G}_W$. One can follow the proof
of Russo's formula in e.g.\ \cite{MM} to obtain for any $g_W \in \mathcal{P}_W$ that
\begin{displaymath}
\frac{d}{dr}\mathbb{P} _{p,r}(A|\mathcal{G}_W=g_W)=\mathbb{E}_{p,r}\left( n(A) |\mathcal{G}_W=g_W\right).
\end{displaymath}
Since
\begin{displaymath}
\frac{d}{dr}\mathbb{P} _{p,r}(A)=\frac{d}{dr}\sum \limits _{g_W\in \mathcal{P}_W} \mathbb{P}
_{p,r}(A|\mathcal{G}_W=g_W)\nu _p(\mathcal{G}_W=g_W).
\end{displaymath}
and
the sum is finite,
the sum and the derivative can be interchanged, giving
\begin{eqnarray*}
 \frac{d}{dr}\mathbb{P} _{p,r}(A)&=& \sum \limits _{g_W\in \mathcal{P}_W} \nu _p(\mathcal{G}_W=g_W)
\frac{d}{dr}\mathbb{P} _{p,r}(A|\mathcal{G}_W=g_W)\\
&=& \sum \limits _{g_W\in \mathcal{P}_W} \nu _p(\mathcal{G}_W=g_W) \mathbb{E} _{p,r}(n(A)|\mathcal{G}_W=g_W)
\\
& =&\mathbb{E}_{p,r}(n(A))
\end{eqnarray*}
\hfill $\Box $\\

\begin{Corollary}
\label{Russocons}
If $A$ is a decreasing event depending on colours of vertices in a finite set $W\subset \mathbb{Z}^2$,
then we have, for all $p\in [0,1]$,
\begin{displaymath}
\frac{d}{dr}\mathbb{P}_{p,r}(A)=-\mathbb{E}_{p,r}(n(A)).
\end{displaymath}
\end{Corollary}

The following lemma gives a finite size criterion for percolation
(see \cite{Russocrpr}, Lemma 2).

\begin{Lemma}\label{finitesizecrit}
There exists a constant $\varepsilon >0$ which satisfies the
following property. If there exists $N\in \mathbb{N}$ such that
\begin{equation}\label{littlerange}
(N+1)(3N+1)\nu _p(\mathcal{D}(0) \geq \frac{N}{3}) \leq \varepsilon
\end{equation}
and
$$
\mathbb{P}_{p,r}(V_{N,3N}^b) >  1-\varepsilon ,
$$
then $\Theta (p,r)>0$.
If there exists $N\in \mathbb{N}$ such that (\ref{littlerange}) holds and
$$
\mathbb{P}_{p,r}(V_{N,3N}^{w*})> 1-\varepsilon,
$$
then $\Theta ^*(p,1-r)>0$.
\end{Lemma}

We do not give the proof of this lemma here as it uses a well-known
coupling argument with a 1-dependent bond percolation model on
$\mathbb{Z}^2$ (see, e.g., the proof of Theorem~2.6 in \cite{DaC}
or the proof of Theorem~1 in \cite{Voronoi}). Using Theorem~\ref{lexpdecay},
Lemma~\ref{finitesizecrit}, and standard arguments, we obtain the
following lemma, which relates the occurrence of percolation to the
probability of crossing large rectangles.

\begin{Lemma}
\label{lessoreq}
For $p<1/2$, we have
\begin{enumerate}
\item[(a)] $\limsup \limits _{n\to \infty } \mathbb{P}_{p,r}(V_{n,3n}^b)=1$ if and only if $\Theta
(p,r)>0$.\\
\item[(b)] $\limsup \limits _{n\to \infty } \mathbb{P}_{p,r}(V_{n,3n}^{w*})=1$ if and only if $\Theta
^*(p,1-r)>0$.
\end{enumerate}
\end{Lemma}

In order to state the final result in this section, taken from~\cite{GKR}, we need the following notation.
Let $\mu $ be a probability measure on colour configurations where the vertices of $\mathbb{Z}^2$ are each
declared black or white.
Let us denote the event that the origin is in an infinite black cluster by $0\leftrightarrow \infty $, and
the number of infinite black clusters by $N$.

\begin{Theorem}\emph{(\cite{GKR})}
\label{GanKeR}
Assume that
\begin{enumerate}
\item[(1)] $\mu $ is invariant under horizontal and vertical translations and axis reflections.
\item[(2)] $\mu $ is ergodic (separately) under horizontal and vertical translations.
\item[(3)] For any increasing events $E$ and $F$,
$$ \mu (E\cap F)\geq \mu (E)\mu (F).$$
\item[(4)] $0<\mu (0\leftrightarrow \infty )<1$.
\end{enumerate}
If assumptions (1)-(4) hold, then
$$ \mu (N=1)=1 .$$
Moreover, any finite set of vertices is surrounded by a black circuit with probability 1 and, equivalently,
all white $*$-clusters are finite with probability 1.
\end{Theorem}

\section{Proof of Theorem \ref{duality}}
\label{mainproof}

In this section, we shall prove that for any $p<p_c$ we have $r_c(p)+r_c^*(p)= 1$.
This can be split in two parts.
The first one is an easy consequence of Theorem~\ref{GanKeR}, stated in the previous section.

\begin{Theorem}
\label{groreq}
For $p<p_c$, we have $r_c(p)+r_c^*(p)\geq 1$.
\end{Theorem}

\noindent \textbf{Proof.} We apply Theorem \ref{GanKeR}. Let us fix $p<p_c$, and assume that
$r_c(p)+r_c^*(p)< 1$.
Then, we may choose some $r\in (r_c(p),1-r_c^*(p))$. Since $r>r_c(p)$, we have $\Theta (p,r)>0$.
On the other hand, it is clear that $\Theta(p,r) \leq r$.
This gives $0<\Theta (p,r)<1$, i.e.\ (4) for the measure $\mathbb{P}_{p,r}$.
Condition (3) is provided by Theorem \ref{FKG}, (2) by Corollary \ref{erg}, and
(1) clearly holds for $\mathbb{P}_{p,r}$. Therefore, all white $*$-clusters
are finite with probability 1. However, this cannot be the case since $r<1-r_c^*(p)$.
\hfill $\Box $\\

To prove the difficult direction, $r_c(p)+r_c^*(p)\leq 1$, we shall use ideas described in~\cite{Russocrpr},
some of which are based on Kesten's proof of $p_c=1/2$ for Bernoulli bond percolation on
$\mathbb{Z}^2$ (see \cite{Kesten}). However, here the proof is considerably more difficult
due to the dependence structure of the DaC model.
Some difficulties are of a geometrical nature, others arise from the fact that we have
to use an RSW type theorem which is weaker than the RSW theorem available for Bernoulli
(independent) percolation, and used by Kesten~\cite{Kesten} in his celebrated proof.

\begin{Theorem}
\label{main}
For any $p<p_c$, $r_c(p)+r_c^*(p)\leq 1$.
\end{Theorem}

\noindent \textbf{Proof.} We shall prove this theorem by contradiction.
Assume that for some $p<p_c$, $r_c(p)+r_c^*(p)> 1$, and fix such a $p$.
Most of the time in the rest of the proof, this $p$ will not appear in our notation.
By the assumption above, we can choose $1-r_c^* < r_1 < r_2 < r_c$. Since for all
$r \in (1-r^*_c,r_c)$, $\Theta (p,r)=0$ and $\Theta ^*(p,1-r)=0$, we have by
Lemma~\ref{lessoreq} that
\begin{displaymath}
\limsup \limits _{n\to \infty } \mathbb{P}_{r_2}(V_{n,3n}^b)<1,
\end{displaymath}
and
\begin{displaymath}
\limsup \limits _{n\to \infty } \mathbb{P}_{r_1}(V_{n,3n}^{w*})<1.
\end{displaymath}
Applying Lemma~\ref{selfdl}, we obtain from these inequalities that
\begin{equation}\label{white}
\liminf \limits _{n\to \infty } \mathbb{P}_{r_2}(H_{n,3n}^{w*})>0,
\end{equation}
and
\begin{equation}\label{black}
\liminf \limits _{n\to \infty } \mathbb{P}_{r_1}(H_{n,3n}^{b})>0.
\end{equation}

Inequality (\ref{white}) implies that there exists $\gamma >0$ and an integer $N_0$ such that for
all $n>N_0$, we have $\mathbb{P}_{r_2}(H_{n,3n}^{w*})\geq \gamma$.
Since $H_{n,3n}^{w*}$ is a decreasing event, by monotonicity this inequality holds in the whole interval:
for all $r\in [r_1,r_2]$ and all $n>N_0$,
\begin{equation}
\label{ivwcr}
\mathbb{P}_{r}\left( H_{n,3n}^{w*}\right )\geq \gamma .
\end{equation}

Since the measure $\mathbb{P}_{r_1}$ is invariant under 90 degree rotations, Theorem \ref{RSW} and
inequality
(\ref{black}) imply that
\begin{equation}\label{black2}
\limsup \limits _{n\to \infty } \mathbb{P}_{r_1}(V_{n,3n}^{b})>0,
\end{equation}
which implies that there exists $\alpha >0$ and a sequence of side lengths $n_k\to \infty $ as
$k\to \infty$ such that for every $k$,
\begin{equation}\label{Mfolott}
\mathbb{P}_{r_1}(V_{n_k,3n_k}^b)\geq \alpha.
\end{equation}

For later purposes we remark that the FKG inequality (Theorem~\ref{FKG}) and a standard pasting
argument imply that for each $k$ we have, for all $i \in \mathbb{N}, i \geq 1$,
\begin{equation}\label{hosszabbak}
\mathbb{P}_{r_1}(V_{n_k,(2i+1)n_k}^b)\geq \alpha ^{2i-1}.
\end{equation}
Indeed, consider rectangles $[0,n_k]\times [2jn_k,(2j+3)n_k]$ for $j=0,1,\ldots ,i-1$ and squares
$[0,n_k]\times [2ln_k,(2l+1)n_k]$ for $l=1,2,\ldots i-1$.
If there are black vertical crossings in these $i$ rectangles of size $n_k\times 3n_k$ and black horizontal
crossings in the $i-1$ squares of size $n_k\times n_k$, then there is a black vertical crossing in the
$n_k\times (2i+1)n_k$ rectangle since horizontal and vertical crossings of the same square meet.
Using Theorem~\ref{FKG} and the fact that the probability of a horizontal black crossing in a square
is bounded below by the probability of a vertical crossing in an $n_k\times 3n_k$ rectangle,
we obtain~(\ref{hosszabbak}).

Let us now fix an integer $L$ with the property that if we consider $L$ Bernoulli (i.e. independent)
trials, each with success probability $\alpha ^4/2$, then the probability that there are at least
$\frac{32}{(r_2-r_1) \gamma \alpha ^{63}}+1$ successes is at least $1/2$.

Next, we choose an element $m_1$ of the sequence $\{ n_k\}$ (for which~(\ref{Mfolott})
holds) so large that it satisfies
\begin{equation}\label{sokkicsi}
(8m_1^2+8m_1)e^{-\psi (p)\frac{m_1}{100}}\leq \min \left( \frac{\alpha ^4}{2L},\frac{\alpha ^{63}}{4(\alpha
^{63}+1)}\right)
\end{equation}
and
\begin{equation}\label{m1nagy}
m_1>\frac{600}{\psi (p)},
\end{equation}
where $\psi (p)$ is the constant corresponding to our fixed $p$ in Theorem \ref{lexpdecay}.
Then take other elements $m_2,m_3,\ldots ,m_L$ in the sequence $\{ n_k\} $ satisfying
\begin{equation}\label{midiff}
m_{i+1}>\frac{151}{49}m_i
\end{equation}
for $i=1,2,\hdots ,L-1$.
Finally, using the constant $N_0$ from (\ref{ivwcr}), we set $N=8n_j$ for some $j$ such that
\begin{equation} \label{Nnagy}
N>\max{(N_0,\frac{9}{2}m_L)}.
\end{equation}
As $N>N_0$, for all $r\in \left[ r_1,r_2\right]$, we have
\begin{equation}\label{whcr}
\mathbb{P}_{r}\left( H_{N,3N}^{w*}\right )\geq \gamma .
\end{equation}

Since the annulus $A_n = S_{3n,3n} \setminus (S_{n,n}+(n,n))$ can be split into four overlapping rectangles,
each with sides of length $n$ and $3n$, a standard argument, based on pasting crossings
and the FKG inequality (Theorem~\ref{FKG}) implies that, for $i=1,2,\ldots,L$,
\begin{displaymath}
\mathbb{P}_{r_1}(B_{m_i}^b)\geq \left( \mathbb{P}_{r_1}(V_{m_i,3m_i}^b)\right)^4 .
\end{displaymath}

Since $m_1,\ldots,m_L$, and $\frac{N}{8}$ are elements of the sequence $\{ n_k\} $,
we get by (\ref{Mfolott}) and (\ref{hosszabbak}) that $\mathbb{P}_{r_1}(V_{m_i,3m_i}^b)\geq \alpha $
and $\mathbb{P}_{r_1}(V_{\frac{N}{8},8N}^b)\geq \mathbb{P}_{r_1}(V_{\frac{N}{8},65\frac{N}{8}}^b)
\geq \alpha ^{63}$.
By monotonicity these inequalities hold in the whole interval $[r_1,r_2]$.
Hence, for all $r\in \left[ r_1,r_2 \right]$ and for $i=1,\ldots,L$, we obtain

\begin{equation}\label{anncross}
\mathbb{P}_{r}(B_{m_i}^b)\geq \alpha ^4,
\end{equation}
and
\begin{equation}\label{sqvercr}
\mathbb{P}_{r}(V_{\frac{N}{8},8N}^b)\geq \alpha ^{63}.
\end{equation}

We have now made all the preparation needed for the essential part of the proof.
In the second part, we shall show that there are uniformly many pivotal clusters
for the event $H_{N,8N}^{w*}$, in expectation, in the interval $r\in [r_1,r_2]$.
More precisely, we will show that for all $r\in \left[ r_1,r_2\right] $, we have
\begin{equation}\label{pivotal}
\mathbb{E}_{r}\left( n\left( H_{N,8N}^{w*}\right) \right) > \frac{1}{r_2-r_1},
\end{equation}
where $n\left( H_{N,8N}^{w*}\right)$ denotes the number of $p$-clusters that are
pivotal for the event $H_{N,8N}^{w*}$.

Before giving the proof, let us explain how this statement leads to a contradiction.
By putting Corollary \ref{Russocons} and (\ref{pivotal}) together, we obtain
\begin{displaymath}
\max \limits _{r\in \left[ r_1,r_2 \right]}\frac{d}{dr}\mathbb{P}_r(H_{N,8N}^{w*})<-\frac{1}{r_2-r_1}.
\end{displaymath}
However, this cannot be the case since it would imply
\begin{eqnarray}
\mathbb{P}_{r_2 }\left( H_{N,8N}^{w*}\right) & \leq & \mathbb{P}_{r_1}\left( H_{N,8N}^{w*}\right) \
 +(r_2-r_1) \max \limits _{r\in \left[ r_1,r_2 \right]}\frac{d}{dr}\mathbb{P}_r(H_{N,8N}^{w*}) \nonumber \\
 & < & \mathbb{P}_{r_1}\left( H_{N,8N}^{w*}\right) - 1, \nonumber
\end{eqnarray}
which is clearly impossible.

Note that it was the assumption that the interval $(1-r_c^*,r_c)$ is non-empty that enabled us to choose
a sub-interval $[r_1,r_2]$ of positive length where the derivative of $\mathbb{P}_r(H_{N,8N}^{w*})$ is
uniformly bounded away from $0$ by $-\frac{1}{r_2-r_1}$.
Since this leads to a contradiction, we conclude that $r_c \leq 1-r_c^*$, as stated in Theorem \ref{main}.
It remains to prove (\ref{pivotal}).\\

\noindent \textbf{Proof of inequality (\ref{pivotal}).}
It will be convenient to introduce the following notation to denote certain parts of $S:=S_{N,8N}$.
We will denote by $t(S):=[0,N]\times \{ 8N \}$ the \textit{top of $S$}, by $b(S):=[0,N]\times \{ 0 \}$
its \textit{bottom}, by $l(S):=\{ 0 \} \times [0,8N]$ its \textit{left side}, by
$r(S):=\{ N \} \times [0,8N]$ its \textit{right side}, by $UH(S):=[0,N]\times [4N+1,8N]$ its
\textit{upper half}, and by $LH(S):=[0,\frac{N}{2}]\times [0,8N]$ its \textit{left half}.

We shall now present a construction of black and white paths in $S$ which guarantees the existence of
many pivotal clusters, and which succeeds with a high enough probability to provide the desired lower
bound for $\mathbb{E}_{r}\left( n\left( H_{N,8N}^w\right) \right)$.
The construction consists of three parts. In the first part, we show that with probability bounded away
from 0, there is a horizontal white $*$-crossing in the lowest part of $S$.\\

\noindent {\bf Part 1.}
We start looking for the lowest white horizontal $*$-crossing of $S$.
It is well-known that the lowest such $*$-crossing can be found (when it exists) by checking only
the colours of vertices (in $S$) below the $*$-crossing and on it.
(The meaning of expressions such as ``below, above, to the right of'' can be made precise via the
Jordan Curve Theorem.)

Recall that by (\ref{whcr}), the probability of the event $H^{w*}_{N,3N}$ is uniformly bounded
below by $\gamma$.
Suppose that $H^{w*}_{N,3N}$ occurs.
Denote the lowest horizontal white $*$-crossing in $S_{N,3N}$ by $\Pi _h^{w*}$.
We shall later use the fact that $\Pi _h^{w*}$ is also the lowest horizontal white $*$-crossing in $S$.
So far we have checked sites only below or on $\Pi _h^{w*}$, but not above it.
However, since the model is dependent, we do have some information above $\Pi _h^{w*}$;
for example that the $p$-clusters of the vertices in $\Pi _h^{w*}$ are white.
Therefore, let us consider the thickened $*$-crossing
\begin{displaymath}
\Gamma _h ^{w*}:=\bigcup _{x\in \Pi _h^{w*}}C_x ^p.
\end{displaymath}

\begin{figure}[h]
\centering
\includegraphics[scale=0.34, trim= 0mm 0mm 0mm 0mm, clip ]{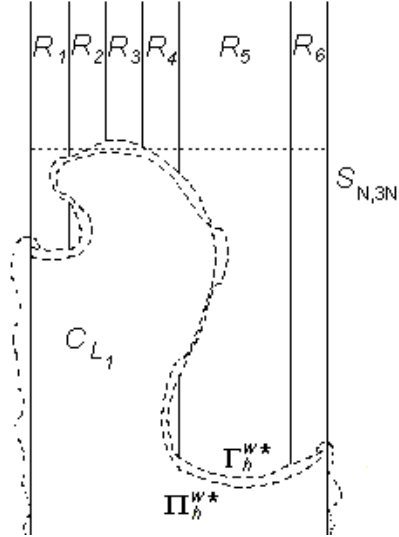}
\caption{Middle part of the rectangle $S$. The top side of $S_{N,3N}$ is
indicated by a dotted segment. The lower broken line represents the lowest
horizontal white $*$-crossing $\Pi^{w*}_h$ of $S_{N,3N}$. The higher broken
line represents the upper boundary of the thickened $*$-crossing $\Gamma^{w*}_h$.
The dotted lines at the left and right side of $S$ form part of the boundary
of $C_{L_1}$.}
\label{rek}
\end{figure}

We denote the portion of $S$ above $\Gamma _h^{w*}$ by $A(\Gamma _h^{w*})$.
We also need to define the following sets (see Figure \ref{rek}):
\begin{eqnarray*}
R_i & := &\left( \left[ (i-1)\frac{N}{8},i\frac{N}{8} \right] \times [0,8N]\right) \cap A(\Gamma _h^{w*})
\textrm{ for } i=1,2,3,4,\\
R_5 & := &\left( \left[ 4\frac{N}{8},7\frac{N}{8} \right] \times [0,8N]\right) \cap A(\Gamma _h^{w*}),\\
R_6 & := &\left( \left[ 7\frac{N}{8},N \right] \times [0,8N]\right) \cap A(\Gamma _h^{w*}),\\
C_{L_1} & := &\bigcup _{x\in S\setminus A(\Gamma _h^{w*})}C_x ^p.
\end{eqnarray*}

We know that the edges in the edge boundary $\Delta _1:=\Delta C_{L_1}$ are closed in the underlying
Bernoulli percolation, hence $\Delta _1$ forms a closed barrier.
Note that this barrier is obtained without checking the states of edges or colours of vertices in
$ext(\Delta _1)$.

We now claim that with high probability, $\{\Gamma _h ^{w*}\cap UH(S)=\emptyset \}$, and
$\{C_{L_1} \cap R_i=\emptyset \}$ for $i=2,3,4,5$.
Indeed, if all the vertices in $S_{N,3N}$ have a dependence range smaller than $N$,
then the first equality holds.
The only way that any of the latter equalities could be false is that there is a vertex
below $\Gamma _h ^{w*}$ whose $p$-cluster extends above $\Gamma _h ^{w*}$ so much that
it intersects one of the rectangles $R_2,R_3,R_4,R_5$.
For this to happen, there has to be a vertex in $S_{N,3N}$ with a dependence range of at
least $\frac{N}{8}$.
Hence, using crude estimations, we give an upper bound for the probability
that at least one equation is false:
\begin{eqnarray*}
\nu _{p}\left( \bigcup_{x\in S_{N,3N}} \{ \mathcal{D}\left( x \right) \geq \frac{N}{8} \} \right) &
\leq & \nu _{p}\left( \bigcup _{x\in S_{N,3N}} \{ \mathcal{D}\left( x \right) \geq \frac{N}{100} \} \right)
\\
& \leq &\sum \limits _{x\in S_{N,3N}}\nu _{p}\left( \mathcal{D}\left( x \right) \geq \frac{N}{100}\right) \\
& \leq & (N+1)(3N+1) \nu _{p}\left( \mathcal{D}\left( 0 \right) \geq \frac{N}{100}\right) \\
& \leq & (8N^2+8N) e^{-\psi (p) \frac{N}{100}} \\
& \leq & (8m_1^2+8m_1) e^{-\psi (p) \frac{m_1}{100}} \\
& \leq & \frac{\alpha ^4}{2L},
\end{eqnarray*}
by the choice of $m_1$ (see (\ref{sokkicsi})). We have also used Theorem \ref{lexpdecay} and the monotonicity of
$f(x)=(8x^2+8x)e^{-\psi (p)\frac{x}{100}}$ for $x>\frac{600}{\psi (p)}$ (this is justified
because of (\ref{m1nagy})).
Therefore, with probability at least $1-\frac{\alpha ^4}{2L}$,  we have no information about the bond
configuration or the colour of the vertices in $R_2,R_3,R_4,R_5$, so their union provides an
unexplored region in $A(\Gamma _h^{w*})$, which contains $UH(S)$.

We now condition on the events $\{\Gamma _h ^{w*}\cap UH(S)=\emptyset \}$, $\{C_{L_1} \cap R_i=\emptyset \}$
for $i=2,3,4,5$, and continue with the second part of our construction.\\

\noindent {\bf Part 2.}
In this part, our task is to find the leftmost vertical black path in $R_2\cup R_3\cup R_4$ from $t(S)$
to $\Gamma_h ^{w*}$.
Here and later, if $W$ is a set of white vertices, then by ``a black path to $W$'' we mean a black path
to some vertex at distance 1 from $W$.
Let us consider the following event:
$E_2:=\{$there is a black path from $t(S)$ to $\Gamma_h ^{w*}$ that does not leave $R_2\cup R_3\cup R_4\}$.
Let us denote the $\sigma $-algebra generated by the information we have so far by $\mathcal{F}_c$, and let
us denote the conditional measure $\mathbb{P}_{r}(\cdot |\mathcal{F}_c)$ by $\mathbb{P}_{r}^{(c)}$.
We shall first show that for all $r\in [r_1,r_2]$, we a.s.\ have
\begin{equation}\label{conditioning1}
\mathbb{P}_{r}^{(c)}(E_2)\geq \frac{\alpha ^{63}}{2}.
\end{equation}

Let $\omega ^{(c)}=(\eta^{(c)},\xi^{(c)})$ and $\omega =(\eta,\xi)$ be elements of $\Omega$
drawn according to $\mathbb{P}_{r}^{(c)}$ and $\mathbb{P}_{r}$, respectively.
We shall show that $\omega$ and $\omega^{(c)}$ can be coupled in
such a way that if there is no large $p$-cluster in $\eta$ in $R_3$,
and there is a vertical black path in $\xi$ from $t(S)$ to $\Gamma_h^{w*}$
that does not leave $R_3$, then there is a black path
in $\xi^{(c)}$ that does not leave $R_2 \cup R_3 \cup R_4$.

We first couple $\eta$ and $\eta^{(c)}$ so that they coincide in the exterior of $\Delta_1$.
This is possible because ${\cal F}_c$ contains information only about $\Delta_1$ and $int(\Delta_1)$,
and Bernoulli percolation configurations restricted to disjoint sets are independent.
Note that the $p$-clusters of $\eta$ can extend beyond $\Delta_1$, therefore each $p$-cluster of
$\eta^{(c)}$ contained in $ext(\Delta_1)$ is a subset of a $p$-cluster of $\eta$, but they are
not necessarily the same.

In order to couple $\xi$ and $\xi^{(c)}$, we consider the collection $\mathbb G$ of all
(self-avoiding) paths in $\left[ \frac{N}{4},\frac{3N}{8} \right] \times [0,8N]$ from
$t(R_3):=\left[ \frac{N}{4},\frac{3N}{8} \right] \times \{8N\}$ to $\Gamma_h^{w*}$, and
give them some deterministic order.
We also order the vertices in each path starting from $t(R_3)$ and ending at $\Gamma_h^{w*}$.
We denote the $j$-th vertex of the $i$-th path by $v _i ^j$.

To each vertex $x \in {\mathbb Z}^2$ we assign a vector $(c_1(x),c_2(x))$, where $c_1$ and
$c_2$ can take three values: black, white or undefined.
Let ${\cal C}_1$ and ${\cal C}_2$ be the collections of all values assigned respectively to
$c_1(x)$ and $c_2(x)$ for all $x \in {\mathbb Z}^2$, indexed by ${\mathbb Z}^2$.
We start with $c_1(x)$ being undefined for all $x \in {\mathbb Z}^2$, and $c_2(x)$ being the
colour of $x$ given ${\cal F}_c$, or undefined if the colour is not known (note that, in
particular, $c_2(x)$ is undefined for all $x \in ext(\Delta_1)$).
We will generate two coupled colour configurations, $\xi$ and $\xi^{(c)}$, according
to the correct marginal distributions with the help of the following algorithm.
Note that the values assigned to $c_1(x)$ and $c_2(x)$ will change, at least for
some $x$, during the algorithmic construction.

Let $c$ be an ``auxiliary'' variable that can take the same three values: black, white and
undefined. We also use two index variables: $i$ and $j$.

\begin{enumerate}
\item $i:=1, \, j:=1$.
\item $c:=c_1(v _i ^j)$.
\item
\begin{itemize}
\item If $c=$ black, $j:=j+1$.
\item If $c=$ white, $i:=i+1$ and $j:=1$. Stop if $i>|\mathbb G|$.
\item If $c=$ undefined, with probability $r$, let $c:=$ black, and with probability $1-r$, let $c:=$ white.
Then set $c_1(x):=c$ for all $x\in C^p _{v _i ^j}(\eta )$ (i.e., for all $x$ in the same $\eta$ $p$-cluster
of the current vertex), and $c_2(x):=c$ for all $x\in C^p _{v _i ^j}(\eta ^{(c)})$.
\end{itemize}
\item Stop if ${\cal C}_1$ contains a black path from $t(R_3)$ to $\Gamma_h^{w*}$, otherwise go back to 2.
\end{enumerate}

After the algorithm stops, we let $\xi(x)=c_1(x)$ for all $x$'s such that $c_1(x)$ is not undefined,
and $\xi^{(c)}(x)=c_2(x)$ for all $x$'s such that $c_2(x)$ is not undefined.
Note that, because of the nature of the algorithm, the vertices that have not been assigned a colour
are naturally split into $p$-clusters (e.g., if $c_1(x)$ is undefined, then $c_1(y)$ is undefined
for all $y$ in the $\eta$ $p$-cluster of $x$).
We then assign colour black with probability $r$ and white with probability $1-r$ to the $p$-clusters
in $\eta$ and in $\eta^{(c)}$ that have not yet been assigned a colour, independently of each other.

We now make three important observations.

(1) First of all, it can be easily seen that the configurations $\xi$ and $\xi^{(c)}$ generated in
the way described above are distributed according to the correct distributions, ${\mathbb P}_r$
and ${\mathbb P}_r^{(c)}$ respectively.

(2) Moreover, before the very last step of the algorithmic procedure, whenever $c_1(x)$ is black for $x$
in $ext(\Delta_1)$, $c_2(x)$ is also black for that same $x$.
This follows from the fact that, because of the coupling between $\eta$ and $\eta^{(c)}$,
a difference between the $p$-clusters of $\eta$ and those of $\eta^{(c)}$ encountered during
the algorithmic construction can only arise when a $p$-cluster of $\eta$ ``crosses'' $\Delta_1$.
In that case, the $\eta$ $p$-cluster possibly reaches more vertices in $ext(\Delta_1)$ than the
$\eta^{(c)}$ $p$-cluster.
If such an $\eta$ $p$-cluster is coloured white, it makes ${\cal C}_1$ ``more white'' than ${\cal C}_2$.
If it is coloured black, the algorithm stops because a black path from $t(R_3)$ to $\Gamma^{w*}_h$
has been generated.
Therefore, before the very last step of the algorithmic procedure, for every vertex $x$ in $ext(\Delta_1)$
such that $c_1(x)$ is black, $c_2(x)$ is also black, and for every vertex such that $c_2(x)$ is undefined,
$c_1(x)$ is either undefined or white.
This implies that if the algorithm stops because ${\cal C}_1$ contains a black path
from $t(R_3)$ to
$\Gamma^{w*}_h$, then also ${\cal C}_2$ contains a black path from $t(R_3)$ to $\Gamma^{w*}_h$.

(3) Finally, at the end of the algorithmic construction described above, $c_2(x)$ can
be black only if $x$ is in $R_3$ or belongs to the $\eta^{(c)}$ $p$-cluster of a vertex
in $R_3$.

Now note that, because of the coupling between $\eta$ and $\eta^{(c)}$, if $x \in R_3$ has
dependence range not larger than $N/8$ in $\eta$, the same is true for the range of $x$ in
$\eta^{(c)}$.
Therefore, if no vertex in $R_3$ has a range larger than $N/8$ in $\eta$, when the
algorithm stops because it found a black path in ${\cal C}_1$ from $t(R_3)$ to $\Gamma^{w*}_h$, by the
previous comment and observations (2) and (3) above, there is
a black path in ${\cal C}_2$ from $t(R_3)$ to $\Gamma^{w*}_h$ contained inside
$R_2 \cup R_3 \cup R_4$.

It follows that, setting
$$T_{R_3}:=\bigcap _{v\in R_3} \{ \mathcal{D}(v)<\frac{N}{8} \},$$
we obtain a.s.
$$\mathbb{P}_{r}^{(c)}(E_2)\geq \mathbb{P}_{r}\left( V_{\frac{N}{8},8N}^b|T_{R_3}\right )\nu _p(T_{R_3}).$$
Elementary calculations show that
$$\mathbb{P}_{r}\left( V_{\frac{N}{8},8N}^b|T_{R_3}\right )\geq \mathbb{P}_{r}\left( V_{\frac{N}{8},8N}^b\
\right )-\nu _p(T_{R_3}^c),$$
where, with a similar computation as in Part 1,
\begin{eqnarray}
\nu_p(T_{R_3}^c)&=&\nu _p\left( \bigcup _{v\in R_3} \{ \mathcal{D}(v)\geq \frac{N}{8} \} \right) \nonumber
\\
& \leq & (\frac{N}{8} +1)(8N+1) \nu _{p}\left( \mathcal{D}\left( 0 \right) \geq \frac{N}{8} \right) \nonumber
\\
& \leq & (8m_1^2+8m_1) e^{-\psi (p) \frac{m_1}{100}} \nonumber \\
&\leq &\zeta \label{nobigcl},
\end{eqnarray}
where $\zeta=\frac{\alpha ^{63}}{4(\alpha ^{63}+1)}$, and in the last step we used (\ref{sokkicsi}).
As $\mathbb{P}_{r}\left( V_{\frac{N}{8},8N}^b\right ) \geq \alpha ^{63}$ (see (\ref{sqvercr})),
and $\zeta <\frac{\alpha ^{63}}{2(\alpha ^{63}+1)}$, this gives
$$\mathbb{P}_{r}^{(c)}(E_2)\geq (\alpha ^{63}-\zeta )(1-\zeta )\geq \alpha ^{63}(1-\zeta )-\zeta
\geq \frac{\alpha ^{63}}{2},$$ proving (\ref{conditioning1}).

\begin{figure}[h]
\centering
\includegraphics[scale=0.19, trim= 0mm 0mm 0mm 0mm, clip ]{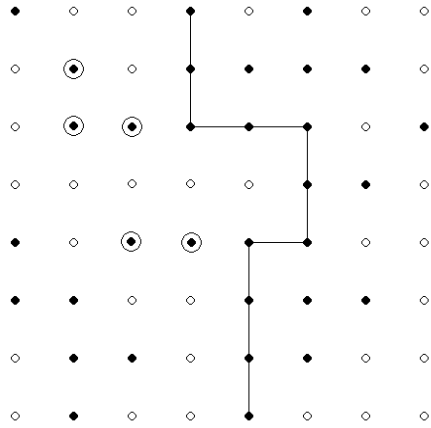}
\caption{The black vertices joined by a line represent a portion of $\tilde\Pi_v^b$.
The black vertices to the left side of $\tilde\Pi_v^b$ and marked by circles belong
to $\Pi_v^b$.}
\label{leftm7}
\end{figure}

Having shown that $\mathbb{P}_{r}^{(c)}(E_2)>0$, conditioning on the event $E_2$, we call
$\tilde\Pi^b_v$ the leftmost black (self-avoiding) path from $t(S)$ to $\Gamma_h^{w*}$
contained in $R_2 \cup R_3 \cup R_4$.
We denote by $\Pi^b_v$ the union of $\tilde\Pi^b_v$ and the black $r$-clusters in $S$ to the left
of $\tilde\Pi^b_v$ connected to $\tilde\Pi^b_v$ (see Figure~\ref{leftm7}).
We denote the region of $S$ to the right of $\tilde\Pi_v^b$ by $J(\tilde\Pi_v^b)$.
Note that no information on the colours of vertices to the right of $\tilde\Pi_v^b$
is needed to determine $\tilde\Pi_v^b$ itself and $\Pi_v^b$.

We now set
$$C_{L_2}:=\bigcup _{x\in (R_2\cup R_3\cup R_4)\setminus (\Pi _v ^b\cup J(\tilde\Pi _v ^b))}C_x^p.$$
If there is no vertex in $R_2\cup R_3\cup R_4$ on the left side of $\tilde\Pi_v^b$
with a dependence range greater than $N$, then $C_{L_2}$ does not intersect $S_{N,7N}\cap J(\tilde\Pi _v ^b)$.
Therefore, the probability of $\{ C_{L_2}\cap S_{N,7N}\cap J(\tilde\Pi_v^b)\neq \emptyset \}$ is at most
$\frac{\alpha ^4}{2L}$, from a computation similar to the one leading to inequality (\ref{nobigcl}).
Conditioning on the event $\{ C_{L_2}\cap S_{N,7N}\cap J(\tilde\Pi _v ^b)= \emptyset \}$, we continue with the
third part of our construction.\\

\noindent{\bf Part 3.} In this part, we shall complete our
construction which shows that there are ``too many'' pivotal
$p$-clusters for having a white horizontal $*$-crossing in $S$. It
is easy to see that if a vertex at distance $1$ from some $x \in
\Gamma^{w*}_h$ is black and is connected to $t(S)$ by a black path
contained in $S$, then $C^p_x$ is pivotal for $H_{N,8N}^{w*}$.
Indeed, due to the fact that every vertex in the lowest horizontal
white $*$-crossing has a black neighbour that is connected to
$b(S)$ by a black path, changing the colour of $C^p_x$ would make
the existence of a white horizontal $*$-crossing of $S$ impossible.

Let $\pi$ be the
``upper layer'' of $\Gamma^{w*}_h$, that is,
the set of vertices in $\Gamma^{w*}_h$ with at least
one neighbour in the region $A(\Gamma^{w*}_h)$.
Denote by $x_0$ the rightmost vertex of $\pi$ at distance 1 from $\Pi _v ^b$.
Furthermore, denote the portion of $\pi$ to the right of $x_0$ by $\Pi_u^{w*}$.
We have specified a site $x_0$ of $\Gamma _h^{w*}$ that lies in $LH(S)\cap S_{N,4N}$
(i.e. the lower-left quarter of $S$) with the property that a black path leads from
$t(S)$ to a point at distance 1 from $x_0$. This implies that $C_{x_0}^p$ is pivotal
for $H_{N,8N}^{w*}$.

Recall that we defined $m_1,\ldots,m_L$ so that
(\ref{sokkicsi})--(\ref{midiff}) and (\ref{anncross}) hold.
Now let us consider annuli $B_i:= \{ c_i+S_{3m_i,3m_i}+\left(
-\frac{3m_i}{2},-\frac{3m_i}{2}\right) \} \setminus \{
c_i+S_{m_i,m_i}+\left( -\frac{m_i}{2},-\frac{m_i}{2}\right) \}$,
$i=1,\ldots,L$, centered at $c_i:=x_0$ if $m_i$ is even and at
$c_i=x_0+(1/2,1/2)$ if $m_i$ is odd. This means that
for any $i$, $B_i$ is an annulus with center at distance at most
$1/2$ from $x_0$, with inner diameter $m_i$ and outer diameter
$3m_i$. We will look for black paths in these annuli between $\Pi
_v^b$ and $\Pi _u ^{w*}$. Note that even the largest annulus, $B_L$
does not go above $S_{N,7N}$, nor to the right of $R_5$, according to
(\ref{Nnagy}). Let us denote the bounded region determined by the
curves $\left( \{ 7\frac{N}{8} \}\times [0,8N]\right) \cap A(\Gamma
_h^{w*}) $ (i.e. the right side of $R_5$), $[0,N]\times \{ 7N \}$
(the top of $S_{N,7N}$), $\Pi _v^b$, and $\Pi _u ^{w*}$, by $AJ$.

Let $\Delta _2$ denote the edge boundary of $C_{L_2}$, defined at the
end of Part 2.
We shall look for black paths in the annuli in the region $AJ$. Let
$B_i[\frac{m_i}{100}]$ denote the $\frac{m_i}{100}$-neighbourhood of
$B_i$. We consider the events $Q_i, i=1,\ldots ,L$, that there is a
black path in $B_i[\frac{m_i}{100}]\cap AJ$ between $\Pi _v^b$ and
$\Pi _u^{w*}$, and $Q:=\{$there are at least $\frac{32}{(r_2-r_1)
\alpha ^{63} \gamma }+1$ indices $i$ such that $Q_i$ holds$\}$. Let
us denote the $\sigma $-algebra generated by the information we have
so far by $\mathcal{F}_{c_2}$, and let us denote the conditional
measure $\mathbb{P}_{r}(\cdot |\mathcal{F}_{c_2})$ by
$\mathbb{P}_{r}^{(c_2)}$. We shall show that for any $r\in
[r_1,r_2]$, a.s.,
\begin{equation}\label{qwerty}
\mathbb{P}_r^{(c_2)}(Q)\geq \frac{1}{4}.
\end{equation}

Let $\omega^{(c_2)}=(\eta^{(c_2)},\xi^{(c_2)})$ and $\omega =(\eta,
\xi)$ be configurations in the plane, drawn according to $\mathbb{P}_{r}^{(c_2)}$
and $\mathbb{P}_{r}$, respectively. We shall show that $\omega $ and
$\omega^{(c_2)}$ can be coupled in such a way that if in $\omega$, for
all $i=1,\ldots,L$, there is no vertex in $B_i$ with a dependence range
larger than $\frac{m_i}{100}$, and for some $j$ there is a black
circuit in $B_j$ in $\xi $, then $\xi ^{(c_2)}\in Q_j$.

First, we couple the edge configurations $\eta $ and $\eta ^{(c_2)}$
in $ext(\Delta _1)\cap ext(\Delta _2)$.
$\mathcal{F}_{c_2}$ contains information that the edges in $\Delta _1\cup\Delta _2$
are closed, and $\Pi _v ^b$ is coloured black
(plus some information about $int(\Delta _1)\cup int(\Delta _2)$,
but that has no influence on the exterior).
This implies that, according to Lemma \ref{stochdom}, there will be a
bias in the configuration $\eta^{(c_2)}$ towards more open edges.
In fact, according to Lemma \ref{stochdom}, in $ext(\Delta _1)\cap ext(\Delta _2)$,
$\eta$ and $\eta^{(c_2)}$ can be coupled so that any closed edge in the latter is
also closed in the former. We pick two such coupled configurations, and concentrating
on $\eta^{(c_2)}$ first, we denote by $H$ the union of $\Pi _v^b$ with the set of vertices
that are connected to $\Pi _v^b$ by an open path in $\eta ^{(c_2)}$. Then $\Delta H$ is a
closed barrier in $\eta ^{(c_2)}$, and by the coupling, this barrier is also closed in $\eta$.
As a final ingredient in our joint construction, we now redraw the configurations in
$ext(\Delta _1)\cap ext(\Delta _2)\cap ext(\Delta H)$ in both configurations, so that in
this region the configurations agree, and are (conditionally) independent of
$\Delta _1, \Delta _2,\Delta H$ and their interiors. The configurations chosen this
way, denoted again by $\eta$ and $\eta^{(c_2)}$, have the correct marginal distributions.

As in Part 2, we assign a vector $(c_1(x),c_2(x))$ to each vertex $x\in \mathbb{Z}^2$,
and let ${\cal C}_1,{\cal C}_2$ be the collections of all the corresponding values,
indexed by ${\mathbb Z}^2$. We take $c_2(x):=$ black for all $x\in int(\Delta H)$.

We shall perform $L$ algorithms, where the $k$-th algorithm corresponds to searching for black paths in $B_k$.
We start with $k=1$.
Let ${\mathbb G}_k$ be the collection of all the (self-avoiding) paths in $B_k\cap AJ$, leading from $\Pi _v
^b$ to $\Pi _u^{w*}$.
We equip ${\mathbb G}_k$ with an arbitrary deterministic ordering. We also order the vertices along each
path, starting from
$\Pi _v ^b$, going towards $\Pi _u^{w*}$. As before, the $j$-th vertex in the $i$-th path in $B_k\cap AJ$ is
denoted by $_k v _i ^j$.
The algorithm which generates $\xi $ and $\xi ^{(c_2)}$ is the same as in Part 2, using the ``auxiliary''
variable $c$ and the index variables $i$ and $j$.

\begin{enumerate}
\item $i:=1,j:=1$.
\item $c:=c_1(_k v _i ^j)$.
\item
\begin{itemize}
\item If $c=$ black, $j:=j+1$.
\item If $c=$ white, $i:=i+1$ and $j:=1$. Stop if $i>|{\mathbb G}_k|$.
\item If $c=$ undefined, with probability $r$, let $c:=$ black, and with probability $1-r$, let
$c:=$ white. Then set $c_1(x):=c$ for all $x\in C^p _{_k v _i ^j}(\eta )$ (i.e., for all $x$ in
the same $\eta$ $p$-cluster of the current vertex), and $c_2(x):=c$ for all
$x\in C^p _{_k v _i ^j}(\eta ^{(c_2)})$.
\end{itemize}
\item Stop if ${\cal C}_1$ contains a black path from $\Pi _v ^b$ to $\Pi _u^{w*}$ in $B_k\cap AJ$,
otherwise go back to 2.
\end{enumerate}

When the algorithm terminates, we increase $k$ by one, and if $k\leq L$, we re-run the
algorithm with the new value of $k$.
After the last algorithm stops, we set $\xi(x)=c_1(x)$ for all $x$ such that $c_1(x)$ is not undefined,
and $\xi^{(c_2)}(x)=c_2(x)$ for all $x$ such that $c_2(x)$ is not undefined.
We then assign colour black with probability $r$ and white with probability $1-r$ to the $p$-clusters
in $\eta$ and in $\eta^{(c)}$ that have not been assigned a colour yet, independently of each other.
Here, we make three important remarks.

(1) Due to the coupling, the bond configurations $\eta $ and $\eta ^{(c_2)}$ are the same
in $ext(\Delta _1)\cap ext(\Delta _2)\cap ext(\Delta H)$.
Therefore,
if for all $k=1,\ldots ,L$, there is no vertex in $B_k\cap AJ\cap ext(\Delta H)$ with
a dependence range larger than $\frac{m_k}{100}$ in $\eta $, then the same is true in $\eta ^{(c_2)}$ as
well.
By re-writing
(\ref{midiff}) as
$\frac{m_{k+1}}{2}>\frac{3m_k}{2}+\frac{m_k}{100}+\frac{m_{k+1}}{100}$,
we see that the
$\frac{m_k}{100}$-neighbourhoods of the annuli $B_k$ are disjoint.
Hence, if for all $k=1,\ldots ,L$, there is no vertex in $B_k$ with
a dependence range larger than $\frac{m_k}{100}$ in $\eta $, then any point gets
$c_1$ or $c_2$ values by at most one of the algorithms.

(2) Similarly to Part 2, the configurations $\xi$ and $\xi^{(c)}$ generated in
the way described above are distributed according to the correct distributions, ${\mathbb P}_r$
and ${\mathbb P}_r^{(c_2)}$ respectively. Note that assigning black to $\xi ^{(c_2)}(x)$ for all
$x\in int(\Delta H)$ is justified since, by the definition of $H$, every such $x$
is connected to $\Pi _v^b$ by an open path in $\eta ^{(c_2)}$.

(3) For any $k$, if there is no vertex in $B_k\cap AJ$ with a dependence range larger than
$\frac{m_k}{100}$, then
before the very last step of the $k$-th algorithmic procedure, whenever $c_1(x)$ is black for $x$
in $ext(\Delta_1)\cap ext(\Delta _2)\cap B_k$, $c_2(x)$ is also black for that same $x$.

To see this, we need to notice that due to the coupling between $\eta$ and $\eta^{(c_2)}$,
the $p$-clusters of a vertex $x\in ext(\Delta_1)\cap ext(\Delta _2)\cap B_k$ in $\eta$ and in $\eta^{(c_2)}$
may differ in the following four cases:
\begin{itemize}
\item $x\in int(\Delta H)$,
\item $C_x ^p (\eta )$ ``crosses'' $\Delta H$,
\item $C_x ^p (\eta )$ ``crosses'' $\Delta_2$,
\item $C_x ^p (\eta )$ ``crosses'' $\Delta_1$.
\end{itemize}
The difference in the first case is unimportant since we have $c_2(x)=$black for all $x\in int(\Delta H)$.
Recall that $\Delta H$ is a closed barrier both in $\eta$ and in $\eta ^{(c_2)}$; hence the
second case never happens. The $k$-th algorithm assigns values to vertices in $B_k\cap AJ$ only.
If there is no vertex in $B_k\cap AJ$ with a dependence range larger than $\frac{m_k}{100}$,
then the third case does not happen either: $\Delta H$ prevents $C_x ^p(\eta )$ for $x\in B_k\cap AJ$
from intersecting $\Delta _2$.
The fourth case is handled exactly the same way as in Part 2: if such a $p$-cluster is coloured white,
it makes ${\cal C}_1$ ``more white'' than ${\cal C}_2$; if it is black, the algorithm has found
an appropriate black path in ${\cal C}_1$ and therefore terminates.

This shows that, for every $k=1,\ldots,L$, if there are no large $\eta $ $p$-clusters in $B_k$,
the presence of a black path in ${\cal C}_1$ from $\Pi _v ^b$ to $\Pi _u^{w*}$ in $B_k\cap AJ$
implies that there is a black path in ${\cal C}_2$ from $\Pi _v ^b$ to $\Pi _u^{w*}$.
Remark (1) above shows that this black path in ${\cal C}_2$ is indeed contained in
$B_k[\frac{m_k}{100}]\cap AJ$.

This implies that if we let $T_{B_i}:=\bigcap _{v\in
B_i}\{\mathcal{D}(v)<\frac{m_i}{100}\},$ $B(B_i):=\{$there is a
black circuit in $B_i$ surrounding $x_0\}$ for $i=1,\ldots ,L$, and
$T:=\bigcap _{i=1} ^LT_{B_i}$, we obtain a.s.
\begin{displaymath}
\mathbb{P}_r^{(c_2)}(Q)\geq \mathbb{P}_{p,r}\ (\textrm{$B(B_i)$
holds for at least $\frac{32}{(r_2-r_1) \alpha ^{63} \gamma }+1$
indices}|T)\nu _p(T).
\end{displaymath}

The second factor is very close to one as
\begin{eqnarray}
1-\nu _p(T) & = & \nu _p(\bigcup _{i=1} ^L T_{B_i}^c) \nonumber \\
& \leq & \sum \limits _{i=1} ^L \sum \limits _{v\in B_i}\nu _p(\mathcal{D}(v)\geq \frac{m_i}{100}) \nonumber \\
& \leq & \sum \limits _{i=1} ^L (8m_i^2+8m_i)\nu _{p}(\mathcal{D}(0)\geq \frac{m_i}{100}) \nonumber \\
& \leq & L\cdot \frac{\alpha ^4}{2L} = \alpha ^4/2 \label{Tbig},
\end{eqnarray}
where we used translation invariance, the monotonicity of the function
$f(x)=(8x^2+8x)e^{-\psi(p)\frac{x}{100}}$ above
$x=\frac{600}{\psi (p)}$, and inequalities (\ref{m1nagy}) and (\ref{sokkicsi}).
Note that, conditioned on $T$, the event $B(B_i)$ depends on the
$\frac{m_i}{100}$-neighbourhood of $B_i$ only. We know the
$\frac{m_i}{100}$-neighbourhoods of the annuli $B_i$ are disjoint.
Therefore, the events $B(B_i)$ ($i=1,\ldots ,L$) are conditionally
independent, conditioned on $T$.
We also have, for $i=1,\ldots ,L$,
$$
\mathbb{P}_{p,r}(B(B_i)|T)\geq \mathbb{P}_{p,r}(B(B_i))-\nu_{p}(T^c)\geq \alpha ^4-\alpha ^4/2,
$$
due to (\ref{anncross}) and (\ref{Tbig}). Hence, by the choice of
$L$ before inequality (\ref{sokkicsi}),
\begin{displaymath}
\mathbb{P}_{p,r}\ (\textrm{$B(B_i)$ holds for at least
$\frac{32}{(r_2-r_1) \alpha ^{63} \gamma }+1$ indices}|T)\geq
1/2.
\end{displaymath}
This shows that
$$
\mathbb{P}_r^{(c_2)}(Q)\geq (1-\frac{\alpha ^4}{2})\frac{1}{2}\geq \frac{1}{4},
$$
proving (\ref{qwerty}).

Note that whenever $Q_i$ happens, there is a pivotal (for the event $H_{N,8N}^{w*}$)
$p$-cluster in or close to $B_i$. Moreover, for any $i\neq j$, the events $Q_i$ and $Q_j$
give rise to different pivotal clusters. Therefore, conditioning on having reached Part 3,
the conditional probability of the event $E_3:=\{$there are at least
$\frac{32}{(r_2-r_1) \alpha ^{63}\gamma }+1$ pivotal clusters for
$H_{N,8N}^{w*}\}$ is at least the conditional probability of $Q$,
which is at least $\frac{1}{4}$, as we have just concluded.

Since it is easy to see that for any $r\in \left[ r_1,r_2 \right]$ the $\mathbb{P}_{r}$-probability of
reaching Part 3 is at least
$\gamma \cdot (1-\frac{\alpha ^4}{2L})\cdot \frac{\alpha ^{63}}{2}\cdot (1-\frac{\alpha ^4}{2L})\geq \frac{\gamma \alpha
^{63}}{8}$, and we know that
$n\left( H_{N,8N}^{w*}\right)$ is a nonnegative random variable, we have for any $r\in \left[ r_1,r_2
\right]$,

\begin{displaymath}
\mathbb{E}_{r}\left( n\left( H_{N,8N}^{w*}\right) \right) \geq \
\left( \frac{32}{(r_2-r_1) \alpha ^{63} \gamma }+1\right) \cdot \frac{1}{4} \cdot \frac{\gamma \alpha
^{63}}{8} >\frac{1}{(r_2-r_1) },
\end{displaymath}
finishing the proof of inequality (\ref{pivotal}), and completing the proof of Theorem \ref{main}.
\hfill $\Box $

\section{Proofs of the remaining results}
\label{consequences}

For the proof of Theorem \ref{criticality}, we need the following result of Russo~\cite{russo1}.
Let $\mu $ be a probability measure that assigns colours black or white to the
vertices of $\mathbb{Z}^2$. Let $P^b _{\infty }(\mu )$ (resp. $P^{w*} _{\infty }(\mu )$) denote the
probability that the black
cluster (resp. white $*$-cluster) of the origin is infinite.
Let $S^b(\mu )$ denote the mean size of the black cluster of the origin.

\begin{Theorem}\emph{(\cite{russo1})}
\label{meansize}
If $\mu $ is translation invariant and $S^b(\mu )<\infty $, then $P^{w*} _{\infty
}(\mu )>0$.
\end{Theorem}

\noindent \textbf{Proof of Theorem \ref{criticality}.} First, we
shall prove criticality when $p<1/2$, $r=r_c(p)$. Our argument
follows the proof of Proposition 1 in \cite{Russocrpr}. Fix $p<1/2$.
Take $\varepsilon >0$ as in Lemma~\ref{finitesizecrit}.
By Lemma~\ref{lexpdecay} and the monotonicity of the function
$f(x)=(x+1)(3x+1)e^{-\psi (p)\frac{x}{3}}$ for $x$ large enough,
there exists $N_0=N_0(p)$ such that, for all $n\geq N_0$,
\begin{equation}\label{nsmalldep}
(n+1)(3n+1)\nu _p(\mathcal{D}(0) \geq \frac{n}{3}) \leq \varepsilon.
\end{equation}
Since $\Theta (p,r)=0$ for all $r<r_c(p)$, Lemma~\ref{finitesizecrit}
and (\ref{nsmalldep}) imply that
$$\mathbb{P}_{p,r}(V_{n,3n} ^b)\leq 1-\varepsilon$$
for all $r<r_c(p),n\geq N_0$.

We claim that for any $n$, the function $\mathbb{P}_{p,r}(V_{n,3n}
^b)$ is continuous in $r$. To see this, notice that the occurrence
of $V_{n,3n}^b$ is completely determined by a partitioning of the
vertices in $S_{n,3n}$ in $p$-clusters, and the colours assigned to
these clusters. Let us denote by $\mathcal{P}_{S}$
the set of partitions of the
vertices in $S_{n,3n}$ which are compatible with a bond configuration,
and the (random)
partition determined by the initial bond percolation by
$\mathcal{G}_{S}$. Fix an arbitrary partition $g_{S}\in
\mathcal{P}_{S}$. Since the colours are assigned independently to
the $p$-clusters determined by $g_{S}$, it is easy to see that
$\mathbb{P}_{p,r}(V_{n,3n} ^b|\mathcal{G}_{S}=g_{S})$ is a
polynomial function of $r$, hence continuous in $r$.
This implies that the (finite) linear combination
$$\mathbb{P}_{p,r}(V_{n,3n} ^b)=\sum \limits _{g_{S}\in \mathcal{P}_{S}}
\nu _p(\mathcal{G}_{S}=g_{S})\ \mathbb{P}_{p,r}(V_{n,3n}
^b|\mathcal{G}_{S}=g_{S})$$ is indeed continuous in $r$.

This shows that for any $n\geq N_0$, if we let $r\to r_c$, we obtain
$\mathbb{P}_{p,r_c}(V_{n,3n} ^b)\leq 1-\varepsilon $. Therefore,
$$\limsup \limits _{n\to \infty }\mathbb{P}_{p,r_c}(V_{n,3n} ^b)<1,$$
which, by Lemma \ref{lessoreq}, implies $\Theta (p,r_c)=0$,
providing the first condition of criticality.

The relation $\Theta ^*(p,r_c^*)=0$ can be proved analogously.
Hence, as $r_c^*=1-r_c,$ we obtain that the
$\mathbb{P}_{p,r_c}$-probability of the origin being in an infinite
white $*$-cluster is 0. Therefore, applying Theorem \ref{meansize}
to the measure $\mathbb{P}_{p,r_c}$ yields that the mean size of
the black cluster of the origin is infinite, concluding the proof
of criticality for $p<1/2$, $r=r_c(p)$.

The fact that there is no infinite black cluster or white
$*$-cluster at $p=1/2$, $r\in (0,1)$,
is a straightforward
consequence of the fact that $\nu_{1/2}$-almost every
$p$-configuration contains infinitely many disjoint open
circuits surrounding the origin. These circuits are coloured
independently, preventing the possibility of black percolation
or white $*$-percolation.
(This idea has been described in \cite{DaC} already
to show that there is no percolation of either colour at $p=1/2, r=1/2$.)
The infinite mean cluster size follows
then from Theorem~\ref{meansize}, as before.

The supercritical case $p>1/2$ is obvious: the probability that the
origin is in an infinite $p$-cluster is positive in that case, and
so is the probability that the colour assigned to that cluster is
black for any $r>0$. \hfill $\Box$\\

\begin{Remark}
{\rm The proof of $r_c+r_c^*=1$ for $p<1/2$ uses the FKG inequality, exponential
decay of correlations and duality.
It is easy to see that polynomial decay of correlations of degree strictly
greater than 2 would be enough for the proof.
The fact that there is no infinite black cluster at $p=1/2$, $r\in(0,1)$,
(Theorem \ref{criticality}) even though duality and the FKG inequality hold
in that case, shows that at $p=1/2$, for any $c>2$ and $N\in \mathbb{N}$,
there exists $n>N$ such that
\begin{displaymath}
\nu _{1/2}(\mathcal{D}(0)\geq n)\geq \frac{1}{n^c},
\end{displaymath}
i.e., in critical bond percolation on the square lattice, the probability that
the origin is connected to $\partial B_n$ by an open path is at least $n^{-c}$.
}
\end{Remark}

\noindent \textbf{Proof of Corollary \ref{phase-diagram}.}
To prove Corollary \ref{phase-diagram}, one needs to put together the results in Theorems
\ref{duality}--\ref{criticality}. Strictly speaking,
the following three statements need additional clarification: for $p<1/2$, we have

(i) $r_c(p) \in [1/2, 1)$,

(ii) $\Theta ^*(p,1-r_c(p))=0$ and the mean size of the white $*$-cluster of the origin is infinite, and

(iii) If $r>r_c(p)$, the size of the white $*$-cluster has an exponentially decaying tail.

\medskip
Now $r_c(p) <1$ follows from Theorem 2.6 in
\cite{DaC}.
The other bound $r_c(p) \geq 1/2$ is an easy consequence of $r_c(p) + r_c^*(p) =1$, since
$r_c(p) \geq r_c^*(p)$.
We have seen the first half of (ii), i.e. $\Theta ^*(p,r_c^*(p))=0$, in the proof of Theorem \ref{criticality}.
We also know that $\Theta(p,r_c(p))=0$, which implies,
according to Theorem \ref{meansize}, that the mean size of the white $*$-cluster of the origin is infinite.
Statement (iii) can be proved the same way as Theorem \ref{exponential-decay}.
\hfill $\Box $\\

The proof of Corollary \ref{continuity} uses the methods of Russo \cite{russo1}, and van den Berg and
Keane \cite{RobKeane},
based on the following lemma, which may be interesting in itself.

\begin{Lemma}
At $p<1/2$, $r>r_c(p)$, the number of infinite black clusters is $\mathbb{P}_{p,r}$-a.s.\ equal to 1.
\end{Lemma}

\noindent \textbf{Proof.}
For $r_c<r<1$, similarly to the proof of Theorem \ref{groreq},
conditions (1)-(4) of Theorem \ref{GanKeR} clearly hold for the measure $\mathbb{P}_{p,r}$.
Theorem~\ref{GanKeR} states that under these conditions, the number of infinite black clusters is 1.
The case $r=1$ is obvious.
\hfill $\Box $\\

\noindent \textbf{Proof of Corollary \ref{continuity}.}
We fix $p<1/2$, and write $\Theta(r|\eta)$ for the conditional probability that the
cluster of the origin is infinite, given that the bond configuration is $\eta$. The above
mentioned classical arguments
and Theorem \ref{criticality} give that $\Theta(r|\eta)$ is for almost all $\eta$ a continuous function
in $r$.

Now fix $\varepsilon >0$ and $r \in [0,1]$. For almost every $\eta$, there exists a maximal $\delta(\eta)$ such that
if $|r'-r| \leq \delta(\eta)$, then $|\Theta(r|\eta)-\Theta(r'|\eta)| < \varepsilon$. Now choose $\bar{\delta}>0$
so small that $\nu_p(\eta; \delta(\eta) < \bar{\delta}) < \varepsilon$ and denote the set $\{\eta; \delta(\eta) \geq
\bar{\delta}\}$ by $A$. Since $\Theta(p,r)=\int \Theta(r|\eta)d\nu_p(\eta)$ we then find that for $r'$ such that
$|r-r'| < \bar{\delta}$, we have
\begin{eqnarray*}
|\Theta(p,r)-\Theta(p,r')| & \leq &  \int_A |\Theta(r|\eta)-\Theta(r'|\eta)| d\nu_p(\eta)
+ \\
& & + \int_{A^c} |\Theta(r|\eta)-\Theta(r'|\eta)| d\nu_p(\eta)\\
&\leq& \varepsilon \nu_p(A) + 2\nu_p(A^c)\\
&\leq & \varepsilon + 2\varepsilon,
\end{eqnarray*}
proving the result.
\hfill $\Box $

\section{The DaC model on the triangular lattice $\mathbb{T}$}
\label{triang}
On the square lattice, the relationship $r_c(p)+r_c^*(p)=1$ does not determine
the critical value $r_c(p)$. However, on the triangular lattice, percolation is
self-dual (i.e., $*$-paths are the same as ordinary paths), so that the same
relationship immediately implies $r_c^*(p) =r_c(p)=1/2$.
In this section, we elaborate a bit on the proof of $r_c(p)+r_c^*(p)=1$ for
$p<p_c(\mathbb{T})$ on the triangular lattice. In this case, the version of the
RSW-type theorem of Bollob\'as and Riordan~\cite{Voronoi} suffices, and we do not
need to use the improvement in~\cite{BBB}.

\begin{figure}[h]
\centering
\includegraphics[scale=0.37, trim= 0mm 0mm 0mm 0mm, clip ]{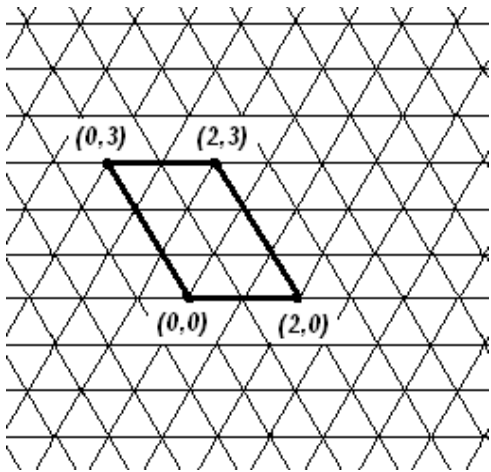}
\caption{Portion of the triangular lattice $\mathbb{T}$. The heavy segments are the sides of
the parallelogram $S_{2,3}$.}
\label{rec}
\end{figure}

We embed the triangular lattice $\mathbb T$ in ${\mathbb R}^2$ so that its vertices $\mathcal{V}(\mathbb T)$
are the intersections of the lines $y=-\sqrt{3}\hspace{0.05cm} x+\sqrt{3}\hspace{0.05cm} k$ and
$y=\frac{\sqrt{3}}{2}\hspace{0.05cm} l$ for $k,l \in \mathbb{Z}$, and denote the elements of
$\mathcal{V}(\mathbb T)$ by $(k,l)$. For example, $\left( 0,0\right)$ refers to the intersection of
$y=-\sqrt{3}\hspace{0.05cm}x$ and $y=0$.
The edges are given by
$\mathcal{E}\left( \mathbb{T} \right) := \{ (a,b): \hspace{0.2cm} a,b\in \mathcal{V}\left( \mathbb{T}\right) ,|\hspace{0.05cm}
a-b\hspace{0.05cm}|=1\} $, where $|\cdot |$ denotes the Euclidean norm (see Figure \ref{rec}).
We define and denote paths, circuits, horizontal and vertical crossings exactly as before.
Note that $S_{m,n}$ corresponds to a parallelogram in $\mathbb{R}^2$ of side lengths $m$ and $n$,
as in the example in Figure~\ref{rec}.

Given the equivalence between crossings and $*$-crossings, we will drop the $*$ from our
notation in this section. We note that the definitions and all the preliminary results of
Section~\ref{prelim} still apply, modulo the reinterpretation of $*$-crossings as ordinary
crossings and the different critical value.
This observation will be implicitly understood in the rest of the section and we will use
the results of Section~\ref{prelim-results} without further comments.
In this section, $p_c$ denotes $p_c(\mathbb{T})$, the critical value for bond percolation on $\mathbb{T}$.

The inequality $r_c(p) \geq 1/2$ can be proved by standard methods.
Similar (however somewhat simpler) considerations to those in the proof of Theorem~\ref{main}
lead to $r_c(p)\leq 1/2$ as follows.
It is easy to see that for any rhombus $S$, the probability of having a black vertical crossing
is exactly the same as the probability of having a black horizontal crossing in $S$.
This observation, together with Lemma~\ref{selfdl} and symmetry of black and white at $r=1/2$
implies the following result.

\begin{Lemma}
\label{rhcross}
For any $p\in [0,1]$ and any $n\in \mathbb{N}$,
\begin{displaymath}
\mathbb{P}_{p,1/2}(V_{n,n}^b)=1/2.
\end{displaymath}
\end{Lemma}

This lemma allows us to use the RSW type theorem of Bollob\'as and Riordan~\cite{Voronoi},
which states that
$\liminf \limits _{n\to \infty }\mathbb{P}_{p,1/2}(V_{n,n}^b)>0$ implies
$\limsup \limits _{n\to \infty }\mathbb{P}_{p,1/2}(V_{n,\rho n}^b)>0$
for any $\rho >0$, to obtain the following result.

\begin{Lemma}
\label{poscross}
For any $p<p_c$, we have
\begin{equation}
\limsup \limits _{n\to \infty }\mathbb{P}_{p,1/2}(V_{n,3n}^b)>0.
\end{equation}
\end{Lemma}

Next we show that certain parallelograms have high crossing probabilities.

\begin{Theorem}
\label{maintri}
For all $p<p_c$, for any $\varepsilon >0$, \hspace{0.05cm}
$$\limsup \limits _{n\to \infty }\mathbb{P}_{p,1/2+\varepsilon }(V_{n,3n}^b)= 1.$$
\end{Theorem}

\noindent \textbf{Proof sketch. }We assume that there exists a $p<p_c$ and an $\varepsilon >0$
such that $\limsup \limits _{n\to \infty }\mathbb{P}_{p,1/2+\varepsilon }(V_{n,3n}^b)<1$.
We denote the measure $\mathbb{P}_{p,r}(\cdot )$ by $\mathbb{P}_{r}(\cdot )$.
The assumption above and Lemma \ref{selfdl} gives $\liminf \limits _{n\to \infty
}\mathbb{P}_{1/2+\varepsilon }(H_{n,3n}^{w})>0$,
which, together with monotonicity, shows that in the interval $r\in [1/2,1/2+\varepsilon ]$,
whenever $n$ is large enough,
the $\mathbb{P}_r$-probability of $H_{n,3n}^{w}$ is bounded away from 0.

On the other hand, Lemma \ref{poscross} gives $\limsup \limits _{n\to \infty
}\mathbb{P}_{1/2}(V_{n,3n}^b)>0$.
Therefore, there exists a sequence of side lengths $n_k\to \infty$ along which
$\mathbb{P}_{1/2}(V_{n_k,3n_k}^b)$ is also bounded away from 0, which, by monotonicity,
gives the same lower bound in the
whole interval $r\in [1/2,1/2+\varepsilon ]$ for $\mathbb{P}_{r}(V_{n_k,3n_k}^b)$.

A careful reading of the proof of Theorem \ref{main} shows that these lower bounds
are enough to determine the parameters $L, m_1,\ldots ,m_L, N$ of the construction described
in the proof of inequality (\ref{pivotal}), which provides a uniform lower bound on the number
of pivotal $p$-clusters for the event $H_{N,8N} ^{w}$ in the interval $r\in [1/2,1/2+\varepsilon ]$,
leading to a contradiction.
\hfill $\Box $\\

Theorem \ref{maintri} together with Lemma \ref{lessoreq} implies
$r_c(p)\leq 1/2$, establishing the equality $r_c(p)=1/2$ and completing the proof of
Theorem~\ref{critical-point}.

We conclude this section with the proof of Proposition~\ref{universality}. \\

\noindent \textbf{Proof of Proposition~\ref{universality}.}
Let us fix $q \geq 2$ and $p<p_c(q)$ and denote the corresponding probability measure
by $\mu_{p,q,r}$.
One can check that conditions (1)--(3) of Theorem~\ref{GanKeR} apply to $\mu_{p,q,1/2}$:
condition (1) is obvious, condition (2) can be found, for example, in~\cite{grimmett2},
condition (3) is proved in~\cite{Haggstrom}.
If we now assume that for $r=1/2$ there exists an infinite black cluster with positive
probability (meaning that condition (4) is also satisfied by $\mu_{p,q,1/2}$), colour
symmetry implies the existence of an infinite white cluster with positive probability,
leading to a contradiction with Theorem~\ref{GanKeR}. We then conclude that there
exists a.s no infinite black cluster at $r=1/2$ and, by colour symmetry again, no
infinite white cluster. Since $\mu_{p,q,1/2}$ is clearly a translation-invariant
measure, Theorem~\ref{meansize} implies infinite mean size for the black $r$-cluster
of the origin.
\hfill $\Box $\\

\medskip\noindent
{\bf Acknowledgements} We would like to thank Rob van den Berg for useful and stimulating discussions.
F.C. thanks Reda J\"urg Messikh and Akira Sakai for interesting discussions at an early stage of this work.


\begin{thebibliography}{99}

\bibitem{ab} M.~Aizenman, D.J.~Barsky, Sharpness of the phase transition in percolation models,
\emph{Comm.~Math.~Phys.} {\bf 86}, 1--48 (1987).

\bibitem{Robnew}
J. van den Berg,
Approximate zero-one laws and sharpness of the percolation transition in a
class of models including 2D Ising percolation,
preprint available from arXiv:math.PR/0707.2077v1 (2007).

\bibitem{BBB}
J. van den Berg, R. Brouwer, B. V\'agv\"olgyi,
Continuity for self-destructive percolation in the plane,
preprint available from arXiv:math.PR/0603223 (2006).

\bibitem{RobKeane}
J. van den Berg, M. Keane,
On the continuity of the percolation probability function, \textit{Particle Systems, Random Media and Large
Deviations} (R.T. Durrett, ed.),
Contemporary Mathematics Series \textbf{26}, AMS, Providence, R. I., 61--65 (1984).

\bibitem{BH}
S.R. Broadbent, J.M. Hammersley, Percolation processes. I. Crystals and mazes, \textit{Proc. Cambridge
Philos. Soc.}
{\bf 53}, 629--641 (1957).

\bibitem{Voronoi}
B. Bollob\'as, O. Riordan,
The critical probability for random Voronoi percolation in the plane is 1/2,
\textit{Probability Theory and Related Fields}
{\bf 136}, 417--468 (2006).

\bibitem{shpthr}
B. Bollob\'as, O. Riordan,
Sharp thresholds and percolation in the plane,
preprint available from arXiv:math.PR/0412510 (2004).

\bibitem{cn1} F.~Camia, C.M.~Newman,
Continuum Nonsimple Loops and 2D Critical Percolation,
\textit{J.~Stat.~Phys.} {\bf 116}, 157-173 (2004).

\bibitem{cn2} F.~Camia, C.M.~Newman,
Two-Dimensional Critical Percolation: the Full Scaling Limit,
\textit{Comm.~Math.~Phys.} {\bf 268}, 1-38 (2006).

\bibitem{cardy} J.L.~Cardy, Critical percolation in finite geometries,
\textit{J.~Phys.~A} {\bf 25}, L201-L206 (1992).

\bibitem{MM}
M. Franceschetti, R. Meester,
\textit{Random networks for communication}, Cambridge University Press (2007).

\bibitem{GKR}
A. Gandolfi, M. Keane, L. Russo,
On the uniqueness of the infinite occupied cluster in dependent two-dimensional site percolation,
\textit{The Annals of Probability} {\bf 16}, 1147--1157 (1988).

\bibitem{Grimmett} G. Grimmett,
\textit{Percolation} (2nd ed.), Springer (1999).

\bibitem{grimmett2} G. Grimmett, \textit{The random-cluster model},
Grundlehren der Mathematischen Wissenschaften [Fundamental Principles of Mathematical Sciences] {\bf 333},
Springer-Verlag, Berlin (2006).

\bibitem{Haggstrom} O. H\"aggstr\"om, Positive correlations in the fuzzy Potts model,
\textit{The Annals of Applied Probability} {\bf 9}, 1149--1159 (1999).

\bibitem{DaC}
O. H\"aggstr\"om,
Coloring percolation clusters at random,
\textit{Stochastic Processes and their Applications} {\bf 96}, 213--242 (2001).

\bibitem{Kesten} H. Kesten, The critical probability of bond percolation on the square lattice equals
$1/2$,
\textit{Comm. Math. Phys.} {\bf 74}, 41--59 (1980).

\bibitem{LigSchSt}
T.M. Liggett, R.H. Schonmann and A.M. Stacey,
Domination by product measures,
\textit{Annals of probability} {\bf 25}, 71--95 (1997).

\bibitem{menshikov} M.V.~Menshikov, Coincidence of critical points in percolation problems,
\emph{Soviet Mathematics Doklady} {\bf 33}, 368--370 (1986).

\bibitem{russo1} L.~Russo, A note on percolation, \textit{Z. Wahrscheinlichkeitstheorie und Verw. Gebiete}
{\bf 43}, 39--48 (1978).

\bibitem{Russocrpr} L. Russo, On the critical percolation probabilities, \textit{Z. Wahrsch. Verw. Gebiete}
{\bf 56}, 229--237 (1981).

\bibitem{smirnov} S.~Smirnov, Critical
percolation in the plane: Conformal invariance,
Cardy's formula, scaling limits,
\textit{C.~R.~Acad.~Sci.~Paris} {\bf 333}, 239--244 (2001).

\bibitem{se1} M.F. Sykes, J.W. Essam, Some exact critical percolation probabilities for site and bond
problems in two dimensions, \textit{Phys. Rev. Lett.} {\bf 10}, 3--4 (1963).

\bibitem{se2} M.F. Sykes, J.W. Essam, Exact critical percolation probabilities for site and bond
problems in two dimensions, \textit{J. Math. Phys.} {\bf 5}, 1117--1127 (1964).

\bibitem{wierman} J.C. Wierman, Bond percolation on honeycomb and triangular lattices,
\textit{Adv. in Appl. Probab.} {\bf 13}, 298--313 (1981).

\end{thebibliography}
\end{document}